\documentclass[11pt]{article}
\usepackage[letterpaper,margin=0.75in]{geometry}
\usepackage{amsmath,amssymb,mathtools}
\usepackage{graphicx}
\usepackage{booktabs}
\usepackage{enumitem}
\usepackage{microtype}
\usepackage{hyperref}
\usepackage{multirow}
\usepackage{longtable}
\usepackage{float}

\newcommand{\erf}{\operatorname{erf}}
\newcommand{\erfc}{\operatorname{erfc}}

\title{\textbf{Structure-Preserving Analysis of a Nonlinear
		Reaction--Diffusion Model for Homogeneous Redox Catalysis}}

\author{
	Guangming Yao$^{a,*}$,
	Michael Shaw$^{b}$,
	Nazia Nimmy$^{a}$,
	Adriana Silva$^{a}$
	\\[0.6em]
	\small
	$^{a}$Department of Mathematics, Clarkson University,
	Potsdam, NY 13699-5815, USA
	\\
	\small
	$^{b}$Department of Chemistry, Southern Illinois University Edwardsville,
	Edwardsville, IL 62026-1652, USA
	\\[0.3em]
	\small
	$^{*}$Corresponding author: gyao@clarkson.edu
}

\date{}

\begin{document}
	\maketitle
	
	\begin{abstract}
		We study a reduced nonlinear reaction--diffusion model for homogeneous redox catalysis under a low-scan-rate, frozen-potential approximation motivated by cyclic voltammetry.
		The reduced system contains two concentration fields coupled through the
		same nonlinear reaction term and therefore satisfies, exactly within the
		reduced model,
		\[
		(a-q)_\tau=(a-q)_{yy}.
		\]
		This identity yields an exact reconstruction of the substrate concentration
		from the electroactive concentration and provides the organizing principle
		for the structure-preserving approximations developed below.

		Formulation I applies the homotopy perturbation method (HPM) to the catalyst field using a constant-coefficient
		Green-function hierarchy and reconstructs the substrate field from the
		exact heat equation. Formulation II linearizes the nonlinear reaction about
		a space- and time-dependent reference profile while preserving the same
		difference structure and the physical far-field state. The methods are
		evaluated using fourth-order method-of-lines discretizations, direct
		Green/heat-kernel formulas for the electrode current, and an independently
		refined nonlinear PDE benchmark. For the baseline problem, the published
		substrate approximation has maximum error \(1.69\times10^{-1}\), whereas
		the \(N=3\) errors are reduced to \(2.09\times10^{-5}\) for Formulation I
		and \(3.61\times10^{-6}\) for Formulation II. At \(N=1\), Formulation II
		improves the errors in the catalyst, substrate, and current by factors of
		approximately \(18.2\), \(46.9\), and \(42.1\), respectively, relative to
		Formulation I. The results show that preserving the exact coupled structure
		substantially improves accuracy and yields an efficient low-order
		analytical framework for nonlinear electrochemical reaction--diffusion
		problems.
	\end{abstract}
	

	\section{Introduction}
	
	Reaction--diffusion models provide a natural mathematical framework for
	electrochemical systems in which interfacial electron transfer, homogeneous
	chemical kinetics, and diffusive mass transport are coupled. Cyclic
	voltammetry has long provided one of the principal tools for probing such
	systems, and its mathematical foundations include the classical treatments
	of reversible, irreversible, and kinetically coupled electrode processes
	\cite{NicholsonShain1964,Nicholson1965,Heinze1984,
		BardFaulkner2001,ComptonBanks2018}. For homogeneous redox catalysis in
	particular, the foundational work of Andrieux and co-workers established the
	general reaction framework and its cyclic-voltammetric consequences
	\cite{AndrieuxPartI1978,AndrieuxEtAl1980}. In these systems, an
	electrochemically generated mediator reacts with a substrate in solution,
	so that the measured current reflects the interaction of electrode kinetics,
	homogeneous reaction kinetics, and transport.
	
	The broader theory of homogeneous molecular electrocatalysis has since been
	developed substantially. General mechanistic and electrochemical treatments
	have clarified how catalytic current, turnover frequency, overpotential,
	intermediate concentrations, and mass-transport effects enter the
	interpretation and benchmarking of homogeneous catalysts
	\cite{Saveant2006,Saveant2008,RountreeEtAl2014,
		CostentinSaveant2014,CostentinSaveant2017,CostentinSaveant2018}.
	Analytical treatments have also been extended to catalytic mechanisms at
	disc electrodes, multielectron and multistep catalysis, reaction-layer
	interpretations, and finite diffusive transport
	\cite{MolinaEtAl2011,LabordaEtAl2022,LabordaEtAl2024,
		MartinezGarciaEtAl2024,Gulaboski2026}. These developments emphasize that
	even relatively compact catalytic mechanisms can lead to nonlinear
	transport problems for which accurate concentration and current expressions
	are nontrivial.
	
	A parallel line of work has focused on analytical and semi-analytical
	solutions of electrochemical reaction--diffusion models. Examples include
	cyclic-voltammetric models for EC-type mechanisms, rotating-disk
	configurations, enzymatic and bioelectrochemical systems, and homogeneous
	catalytic schemes
	\cite{EswariRajendran2011,RajendranKirthigaLaborda2017,
		SaravanakumarEtAl2019,SalomiRajendran2022,
		SylviaRajendran2022,ManimegalaiRajendran2022}.
	More recent steady-state analysis has continued this effort
	\cite{UshaRaniEtAl2023}. A subsequent mathematical comment on that work
	\cite{Tisdell2026} also highlights a point that is especially relevant here:
	closed-form or semi-analytical expressions should be checked not only
	numerically, but also directly against the governing equations and boundary
	conditions.
	
	The homotopy perturbation method (HPM), introduced and developed by He
	\cite{He1999,He2003}, provides one framework for constructing successive
	approximations to nonlinear differential equations. HPM has been applied to
	coupled nonlinear reaction--diffusion systems
	\cite{GanjiSadighi2006}, and convergence of HPM approximations for partial
	differential equations has also been investigated
	\cite{BiazarGhazvini2009}. Within electrochemical modeling, HPM and related
	modified homotopy constructions have been used to derive approximate
	concentration profiles and current responses for enzymatic electrodes,
	convection--diffusion models, biosensors, homogeneous reactions, and cyclic
	voltammetry
	\cite{ShanmugarajanEtAl2011,RasiRajendranSubbiah2015,
		JansiRaniEtAl2017,SwaminathanEtAl2019,SwaminathanEtAl2021,
		EswariRajendran2011,SaravanakumarEtAl2019,SylviaRajendran2022}.
	These studies demonstrate the usefulness of HPM in electrochemical
	reaction--diffusion problems, while also motivating careful examination of
	how the auxiliary homotopy problem is chosen and how the resulting
	approximations respect the structure of the original coupled system.
	
	Accurate numerical simulation remains essential for assessing such
	semi-analytical approximations. Numerical treatments of electrochemical
	reaction--diffusion systems have included collocation methods, stability
	analysis of finite-difference schemes, high-order boundary-gradient
	approximations, and high-order spatial discretizations
	\cite{HertlSpeiser1987,Bieniasz1993,Bieniasz2003,Bieniasz2004}.
	Nonlinear transient reaction--diffusion problems from electroanalytical
	chemistry have likewise been studied through combined analytical and
	numerical approaches \cite{VynnyckyMcKeeBieniasz2021}. These numerical
	considerations are particularly important for electrode-current prediction,
	because the current depends on a concentration derivative at the electrode
	and can therefore be more sensitive to spatial resolution than the
	concentration field itself.
	
	The present work is motivated specifically by the low-scan-rate homogeneous
	redox-catalysis model of Sylvia and Rajendran
	\cite{SylviaRajendran2022}. Their analysis provides closed-form HPM-based
	approximations for the catalyst concentration, substrate concentration, and
	current response. The reduced model contains two concentration fields,
	\(q(y,\tau)\) and \(a(y,\tau)\), coupled through the same nonlinear reaction
	term. This coupling reveals an exact property of the reduced model that is
	independent of any perturbation method: subtracting the two governing
	equations gives
	\[
	(a-q)_\tau=(a-q)_{yy}.
	\]
	Thus the difference between the catalyst and substrate fields satisfies a
	linear heat equation even though the individual fields satisfy nonlinear
	reaction--diffusion equations.
	
	This observation changes the approximation problem in an important way.
	The difference equation is not merely an auxiliary identity; it provides a
	necessary structural consistency condition for any pair of approximate
	concentration fields intended to represent the reduced system. The
	published catalyst approximation provides a useful leading reference
	solution, but its companion substrate approximation does not, in general,
	satisfy this exact difference equation. Consequently, improving the
	approximation is not solely a matter of adding higher-order perturbation
	terms. It is also necessary to preserve the coupling structure already
	present in the governing equations.
	
	The principal contribution of the present work is therefore
	structure-preserving rather than method-specific. We first use the exact
	difference heat equation to derive an explicit reconstruction of the
	substrate concentration from the catalyst concentration. HPM is then used
	as a systematic framework for constructing approximations that respect this
	reconstruction. In Formulation~I, HPM is applied to the catalyst equation
	through a constant-coefficient Green-function hierarchy, while every
	substrate coefficient is recovered from the exact difference relation. We
	show that applying the same frozen-coefficient homotopy directly to both
	governing equations produces the same continuous approximation, providing
	an independent consistency check rather than a distinct method. In
	Formulation~II, the nonlinear reaction is instead linearized about a
	space- and time-dependent reference profile, producing a genuinely
	different correction hierarchy while preserving the same exact difference
	structure and the physical far-field state.
	
	The resulting formulations are assessed against an independently refined
	numerical solution of the original nonlinear reduced system. The numerical
	calculations are used primarily to determine how much accuracy is gained by
	preserving the exact structure, how rapidly the two approximation
	hierarchies converge, and how reliably the electrode current can be
	recovered from the surface concentration gradient. Particular attention is
	given to separating analytical truncation error from spatial
	discretization, quadrature, and reference-solution errors.

	\section{Electrochemical Model and Model Reduction}
	\label{sec:model}
	
	We begin with the homogeneous redox-catalysis mechanism underlying the
	reduced reaction--diffusion system. Homogeneous redox catalysis and its
	cyclic-voltammetric description have been developed extensively, from the
	classical formulations of Andrieux and co-workers to modern treatments of
	molecular electrocatalysis
	\cite{AndrieuxPartI1978,AndrieuxEtAl1980,Saveant2006,
		Saveant2008,RountreeEtAl2014,CostentinSaveant2017}.
	The purpose of this section is to make explicit the connection between this
	electrochemical reaction mechanism and the two-variable dimensionless model
	analyzed in the remainder of the paper.
	
	\subsection{Reaction mechanism}
	
	The electrochemical system follows the homogeneous redox-catalysis
	framework in which an electrochemically generated form of a mediator reacts
	with a substrate in solution
	\cite{AndrieuxPartI1978,AndrieuxEtAl1980,RountreeEtAl2014,
		SylviaRajendran2022}. It consists of an electrode reaction involving the
	substrate couple $A/B$ and a catalytic redox couple $P/Q$. The reaction
	scheme may be written as
	
	\begin{equation*}
		A+e^- \rightleftharpoons B,
	\end{equation*}
	followed by the homogeneous chemical reaction
	\begin{equation*}
		B \xrightarrow{k} C,
	\end{equation*}
	together with the catalytic electrode reaction
	\begin{equation*}
		P+e^- \rightleftharpoons Q,
	\end{equation*}
	and the homogeneous redox-catalytic reaction
	\begin{equation*}
		Q+A
		\mathop{\rightleftharpoons}^{k_1}_{k_2}
		P+B.
	\end{equation*}
	
	Here $P$ and $Q$ denote the oxidized and reduced forms of the catalyst,
	respectively, while $A$ and $B$ denote the corresponding substrate
	species. The active form $Q$ is generated electrochemically at the
	electrode and subsequently reacts with the substrate $A$ in solution.
	The homogeneous reaction regenerates $P$, thereby producing the catalytic
	cycle. This type of mediated catalytic mechanism is a standard setting in
	homogeneous molecular electrocatalysis
	\cite{Saveant2008,RountreeEtAl2014,CostentinSaveant2014,
		CostentinSaveant2018}.
	In the potential region associated with the catalyst wave, the $P/Q$
	electron-transfer reaction is assumed to be sufficiently rapid that the
	surface concentrations satisfy the Nernst relation
	\cite{BardFaulkner2001,ComptonBanks2018}.
	
	\subsection{Dimensional reaction--diffusion equations}
	
	Let
	\[
	c_P(x,t),\qquad
	c_Q(x,t),\qquad
	c_A(x,t),\qquad
	c_B(x,t)
	\]
	denote the dimensional concentrations of the four species, where $x$ is
	the distance from the planar electrode and $t$ is time.
	
	Following the planar diffusion--reaction description used for homogeneous
	redox catalysis and related electrochemical reaction--diffusion models
	\cite{BardFaulkner2001,SylviaRajendran2022,
		VynnyckyMcKeeBieniasz2021}, the model assumes that all species have the
	same diffusion coefficient $D$, diffusion is one-dimensional and
	semi-infinite, and the homogeneous catalytic reaction obeys mass-action
	kinetics. The forward and backward reaction rates are therefore
	\begin{equation*}
		r_f=k_1c_Ac_Q,
		\qquad
		r_b=k_2c_Bc_P.
	\end{equation*}
	
	The net homogeneous catalytic reaction rate is
	\begin{equation*}
		r_{\mathrm{cat}}
		=
		k_1c_Ac_Q-k_2c_Bc_P.
	\end{equation*}
	
	The dimensional reaction--diffusion equations are consequently
	\begin{equation*}
		\frac{\partial c_P}{\partial t}
		=
		D\frac{\partial^2c_P}{\partial x^2}
		+k_1c_Ac_Q-k_2c_Bc_P,
	\end{equation*}
	\begin{equation*}
		\frac{\partial c_Q}{\partial t}
		=
		D\frac{\partial^2c_Q}{\partial x^2}
		-k_1c_Ac_Q+k_2c_Bc_P,
	\end{equation*}
	\begin{equation*}
		\frac{\partial c_A}{\partial t}
		=
		D\frac{\partial^2c_A}{\partial x^2}
		-k_1c_Ac_Q+k_2c_Bc_P,
	\end{equation*}
	and
	\begin{equation*}
		\frac{\partial c_B}{\partial t}
		=
		D\frac{\partial^2c_B}{\partial x^2}
		+k_1c_Ac_Q-k_2c_Bc_P-kc_B.
	\end{equation*}
	
	The first two equations contain equal and opposite homogeneous reaction
	terms because the catalytic reaction converts $Q$ to $P$. Similarly,
	consumption of $A$ is accompanied by formation of $B$. The additional
	term $-kc_B$ represents the irreversible follow-up reaction of $B$.
	
	An important consequence of the first two equations is conservation of
	the total catalyst concentration. Adding the first two equations gives
	\begin{equation*}
		\frac{\partial}{\partial t}(c_P+c_Q)
		=
		D\frac{\partial^2}{\partial x^2}(c_P+c_Q).
	\end{equation*}
	With spatially uniform initial total catalyst concentration and compatible
	boundary conditions, this gives
	\begin{equation*}
		c_P+c_Q=c_P^0,
	\end{equation*}
	where $c_P^0$ is the reference catalyst concentration.
	\subsection{Electrode potential and current}
	
	In cyclic voltammetry, the electrode potential is scanned linearly,
	\begin{equation*}
		E(t)=E_i-vt,
	\end{equation*}
	where \(E_i\) is the initial potential and \(v\) is the scan rate. For a
	Nernstian \(P/Q\) couple, the surface concentrations satisfy
	\begin{equation*}
		\frac{c_Q(0,t)}{c_P(0,t)}
		=
		\exp\left[
		-\frac{F(E(t)-E_{PQ}^{0'})}{RT}
		\right],
	\end{equation*}
	where \(E_{PQ}^{0'}\) is the formal potential of the catalytic couple
	\cite{BardFaulkner2001,ComptonBanks2018}.
	
	Using catalyst conservation at the electrode,
	\[
	c_P(0,t)+c_Q(0,t)=c_P^0,
	\]
	gives
	\begin{equation*}
		\frac{c_Q(0,t)}{c_P^0}
		=
		\frac{1}{1+\exp[-\xi(t)]},
		\qquad
		\xi(t)
		=
		\frac{F(E(t)-E_{PQ}^{0'})}{RT}.
	\end{equation*}
	
	With the dimensionless time
	\[
	\tau=\frac{Fvt}{RT},
	\]
	the linear scan relation implies
	\begin{equation*}
		\xi(\tau)=\xi_i-\tau,
		\qquad
		\xi_i=
		\frac{F(E_i-E_{PQ}^{0'})}{RT}.
	\end{equation*}
	Thus, for a continuously scanned potential, the exact Nernstian boundary
	condition would be time dependent:
	\begin{equation}
		q(0,\tau)
		=
		m(\tau)
		=
		\frac{1}{1+\exp[-(\xi_i-\tau)]}.
		\label{eq:time-dependent-m}
	\end{equation}
	
	The closed-form low-scan-rate approximation studied in the present work,
	however, is based on a constant surface value. Accordingly, for each
	analytical subproblem we regard the dimensionless potential \(\xi\) as a
	fixed parameter and define
	\begin{equation*}
		m_0=m(\xi)
		=
		\frac{1}{1+\exp(-\xi)}.
	\end{equation*}
	The boundary condition used in the structure-preserving analysis is
	therefore
	\begin{equation}
		q(0,\tau)=m_0,
		\label{eq:frozen-boundary}
	\end{equation}
	with \(m_0\) held constant during that transient calculation. We refer to
	this as the \emph{frozen-potential} or quasi-static boundary
	approximation. When dependence on electrode potential is of interest,
	different values of \(\xi\), and hence \(m_0\), may be considered
	parametrically through separate frozen-potential solutions.
	
	It is important to distinguish \eqref{eq:frozen-boundary} from the fully
	time-dependent Nernst boundary \eqref{eq:time-dependent-m}. The analysis
	below concerns the former low-scan-rate reduced problem and does not claim
	to solve the nonautonomous boundary-value problem generated by a continuously
	varying Nernst surface concentration.
	
	The dimensional Faradaic current is determined by the interfacial catalyst
	flux,
	\begin{equation*}
		i
		=
		FSD
		\left.
		\frac{\partial c_P}{\partial x}
		\right|_{x=0},
	\end{equation*}
	where \(S\) is the electrode area
	\cite{BardFaulkner2001,ComptonBanks2018}.
	
	\subsection{Nondimensionalization}
	
	We use the diffusion--voltammetric scaling associated with the reduced
	low-scan-rate model of Sylvia and Rajendran
	\cite{SylviaRajendran2022}. Introduce the dimensionless variables
	\begin{equation*}
		p=\frac{c_P}{c_P^0},
		\qquad
		q=\frac{c_Q}{c_P^0},
		\qquad
		a=\frac{c_A}{c_P^0},
		\qquad
		b=\frac{c_B}{c_P^0},
	\end{equation*}
	together with
	\begin{equation*}
		y
		=
		x\sqrt{\frac{Fv}{RTD}},
		\qquad
		\tau
		=
		\frac{Fvt}{RT}.
	\end{equation*}
	
	The dimensionless bulk substrate concentration is
	\begin{equation*}
		\gamma
		=
		\frac{c_A^0}{c_P^0}.
	\end{equation*}
	
	The dimensionless kinetic parameters are
	\begin{equation*}
		\lambda_1
		=
		\frac{RT}{F}
		\frac{k_1c_P^0}{v},
		\qquad
		\lambda_2
		=
		\frac{RT}{F}
		\frac{k_2c_P^0}{v},
		\qquad
		\lambda
		=
		\frac{RT}{F}
		\frac{kc_P^0}{v}.
	\end{equation*}
	
	With these definitions, the dimensional system becomes
	\begin{align*}
		p_\tau
		&=
		p_{yy}
		+\lambda_1aq
		-\lambda_2bp,
		\\
		q_\tau
		&=
		q_{yy}
		-\lambda_1aq
		+\lambda_2bp,
		\\
		a_\tau
		&=
		a_{yy}
		-\lambda_1aq
		+\lambda_2bp,
		\\
		b_\tau
		&=
		b_{yy}
		+\lambda_1aq
		-\lambda_2bp
		-\lambda b.
	\end{align*}
	
	The catalyst-conservation relation becomes
	\begin{equation*}
		p+q=1.
	\end{equation*}
	\subsection{Reduction to the two-variable model}
	
	We next reduce the four-species transport system to the two concentration
	fields used throughout the subsequent analysis. The catalyst-conservation
	relation gives
	\[
	p=1-q.
	\]
	A second approximation is required to eliminate the intermediate species
	\(b\). Specifically, we impose the local algebraic kinetic closure
	\begin{equation}
		\lambda_1aq-\lambda_2bp-\lambda b=0.
		\label{eq:b-closure}
	\end{equation}
	Equation \eqref{eq:b-closure} assumes that, for purposes of determining
	\(b\) in the reduced model, its local homogeneous production and
	consumption rates are in quasi-steady balance. It is therefore an
	additional modeling approximation and is not an exact consequence of the
	four reaction--diffusion equations above. In particular, it neglects the
	transient and diffusive contributions in the \(b\)-equation when \(b\) is
	eliminated algebraically.
	
	Solving \eqref{eq:b-closure} for \(b\) gives
	\begin{equation*}
		b
		=
		\frac{\lambda_1aq}
		{\lambda+\lambda_2p}.
	\end{equation*}
	Using \(p=1-q\),
	\begin{equation*}
		b
		=
		\frac{\lambda_1aq}
		{\lambda+\lambda_2(1-q)}.
	\end{equation*}
	
	The net reaction term appearing in the \(q\) and \(a\) equations is
	therefore
	\begin{align*}
		\lambda_1aq-\lambda_2bp
		&=
		\lambda_1aq
		-
		\lambda_2(1-q)
		\frac{\lambda_1aq}
		{\lambda+\lambda_2(1-q)}
		\\
		&=
		\lambda_1aq
		\left[
		1-
		\frac{\lambda_2(1-q)}
		{\lambda+\lambda_2(1-q)}
		\right]
		\\
		&=
		\frac{\lambda_1aq}
		{1+(\lambda_2/\lambda)(1-q)}.
	\end{align*}
	
	Defining
	\begin{equation}
		R(a,q)
		=
		\frac{\lambda_1aq}
		{1+\rho(1-q)},
		\qquad
		\rho=\frac{\lambda_2}{\lambda},
		\label{eq:R}
	\end{equation}
	gives the reduced reaction--diffusion system
	\begin{equation}
		\begin{aligned}
			q_\tau &= q_{yy}-R(a,q),\\
			a_\tau &= a_{yy}-R(a,q),
		\end{aligned}
		\qquad y>0,\quad \tau>0.
		\label{eq:system}
	\end{equation}
	
	All structural identities and structure-preserving approximations developed
	below are exact with respect to the reduced two-variable system
	\eqref{eq:system}. They should therefore be distinguished from statements
	about the unreduced four-species transport model, for which
	\eqref{eq:b-closure} is an approximation.
	
	The initial conditions are
	\begin{equation}
		q(y,0)=0,
		\qquad
		a(y,0)=\gamma,
		\label{eq:IC}
	\end{equation}
	and the frozen-potential electrode conditions are
	\begin{equation*}
		q(0,\tau)=m_0,
		\qquad
		a_y(0,\tau)=0.
	\end{equation*}
	The far-field conditions are
	\begin{equation}
		q(y,\tau)\rightarrow0,
		\qquad
		a(y,\tau)\rightarrow\gamma,
		\qquad
		y\rightarrow\infty.
		\label{eq:BCinf}
	\end{equation}
	
	With the sign convention adopted here, the dimensionless current is
	\begin{equation*}
		\psi(\tau)=-q_y(0,\tau).
	\end{equation*}
	Thus the experimentally relevant response depends on a boundary derivative,
	making accurate representation of the near-electrode concentration profile
	particularly important \cite{BardFaulkner2001,Bieniasz2003}.
	
	\section{Exact Structural Reduction}
	\label{sec:structural}
	
	A key feature of \eqref{eq:system} is that the same nonlinear reaction
	term appears in both equations. This produces an exact cancellation that
	is useful both analytically and numerically.
	
	Define
	\begin{equation*}
		w=a-q.
	\end{equation*}
	Subtracting the $q$ equation from the $a$ equation in
	\eqref{eq:system} gives
	\begin{equation}
		w_\tau=w_{yy}.
		\label{eq:wheat}
	\end{equation}
	This relation is exact within the reduced model and does not depend on any
	HPM approximation or parameter limit.
	From \eqref{eq:IC} and \eqref{eq:BCinf},
	\begin{equation*}
		w(y,0)=\gamma,
		\qquad
		w(y,\tau)\rightarrow\gamma
		\quad (y\rightarrow\infty).
	\end{equation*}
	At the electrode,
	\begin{equation*}
		w_y(0,\tau)
		=
		a_y(0,\tau)-q_y(0,\tau)
		=
		-q_y(0,\tau).
	\end{equation*}
	
	Therefore, once $q$ is known, $w$ is determined by a linear heat equation
	and
	\begin{equation*}
		a=q+w.
	\end{equation*}
	
	To obtain an explicit representation, set
	\begin{equation*}
		z=w-\gamma.
	\end{equation*}
	Then
	\begin{equation*}
		z_\tau=z_{yy},
		\qquad
		z(y,0)=0,
		\qquad
		z(\infty,\tau)=0,
		\qquad
		z_y(0,\tau)=-q_y(0,\tau).
	\end{equation*}
	
	Applying the Laplace transform in $\tau$ to this half-line heat problem
	and subsequently inverting the transform gives
	\begin{equation*}
		w(y,\tau)
		=
		\gamma
		+
		\int_0^\tau
		\frac{q_y(0,\sigma)}
		{\sqrt{\pi(\tau-\sigma)}}
		\exp\left[
		-\frac{y^2}{4(\tau-\sigma)}
		\right]\,d\sigma.
	\end{equation*}
	Consequently,
	\begin{equation}
		a(y,\tau)
		=
		q(y,\tau)+\gamma
		+
		\int_0^\tau
		\frac{q_y(0,\sigma)}
		{\sqrt{\pi(\tau-\sigma)}}
		\exp\left[
		-\frac{y^2}{4(\tau-\sigma)}
		\right]\,d\sigma.
		\label{eq:arep}
	\end{equation}
	
	Equation \eqref{eq:arep} provides an exact reconstruction of the substrate
	concentration from $q$. It is central to the first HPM formulation
	developed below: rather than approximating $q$ and $a$ independently, we
	approximate $q$ and recover $a$ in a manner that preserves
	\eqref{eq:wheat} exactly relative to the chosen approximation for $q$.

	\section{Structure-Preserving HPM Formulations}
	\label{sec:hpm}
	
	We first record the published low-scan-rate approximation as a reference
	baseline and then develop two structure-preserving HPM formulations.
	Formulation I uses a constant-coefficient scalar hierarchy for \(q\)
	together with exact structural reconstruction of \(a\). The corresponding
	coupled frozen-coefficient construction is considered separately only to
	establish its equivalence with Formulation I and to provide an independent
	consistency check. Formulation II changes the auxiliary problem by using a
	profile-dependent Taylor linearization of the full nonlinear reaction term.
	
	\subsection{Published Low-Scan-Rate Approximation: Reference Baseline}
	\label{sec:published-baseline}
	
	For comparison with the formulations developed below, we first record the
	published low-scan-rate approximation of Sylvia and Rajendran\,\cite{SylviaRajendran2022}. Writing
	\[
	R(a,q)
	=
	\frac{\lambda_1a}{1+\rho(1-q)}\,q,
	\]
	the published approximation replaces the concentration-dependent
	coefficient multiplying \(q\) by the constant
	\begin{equation*}
		K
		=
		\frac{\gamma\lambda_1}{1+\rho(1-m_0)}.
	\end{equation*}
	This constant coincides with the reference coefficient \(\alpha\)
	introduced in Formulation I.
	The resulting concentration profile is
	\begin{equation*}
		\begin{aligned}
			q_{\mathrm{P}}(y,\tau)
			=
			\frac{m_0}{2}
			\Bigg[
			&
			e^{-y\sqrt{K}}
			\operatorname{erfc}
			\left(
			\frac{y}{2\sqrt{\tau}}-\sqrt{K\tau}
			\right)
			+
			e^{y\sqrt{K}}
			\operatorname{erfc}
			\left(
			\frac{y}{2\sqrt{\tau}}+\sqrt{K\tau}
			\right)
			\Bigg].
		\end{aligned}
	\end{equation*}
	
	This function satisfies the linear diffusion--reaction problem
	\begin{equation*}
		(q_{\mathrm{P}})_\tau
		=
		(q_{\mathrm{P}})_{yy}
		-
		Kq_{\mathrm{P}},
	\end{equation*}
	together with
	\begin{equation*}
		q_{\mathrm{P}}(y,0)=0,
		\qquad
		q_{\mathrm{P}}(0,\tau)=m_0,
		\qquad
		q_{\mathrm{P}}(\infty,\tau)=0.
	\end{equation*}
	
	The published approximation for $a$ is
	\begin{equation*}
		\begin{aligned}
			a_{\mathrm{P}}(y,\tau)
			={}&
			q_{\mathrm{P}}(y,\tau)+\gamma
			-
			m_0
			\operatorname{erfc}
			\left(
			\frac{y}{2\sqrt{\tau}}
			\right)
			\left[
			\sqrt{\pi K\tau}\,
			\operatorname{erf}(\sqrt{K\tau})
			+
			e^{-K\tau}
			\right].
		\end{aligned}
	\end{equation*}
	
	The published substrate expression satisfies the desired electrode condition
	\begin{equation*}
		(a_{\mathrm{P}})_y(0,\tau)=0.
	\end{equation*}
	
	However, in general,
	\begin{equation*}
		(a_{\mathrm P}-q_{\mathrm P})_\tau
		\neq
		(a_{\mathrm P}-q_{\mathrm P})_{yy}.
	\end{equation*}
	Thus the published expression for \(a_{\mathrm P}\) is a separate
	approximation for the substrate concentration rather than the exact
	substrate field implied by \(q_{\mathrm P}\) through the structural
	relation \eqref{eq:arep}.
	
	\subsection{Formulation I: $q$-First Structure-Preserving HPM}
	\label{sec:formulation-I}
	
	Formulation I applies the homotopy perturbation method (HPM) to the \(q\) equation only. The substrate
	concentration is then recovered from the exact structural relation
	\eqref{eq:arep}. This preserves the simplicity of a scalar HPM hierarchy
	while ensuring that every truncated pair satisfies the difference heat
	equation.
	
	\subsubsection{Homotopy construction and reaction coefficients}
	
	Introduce the constant reference coefficient
	\begin{equation*}
		\alpha
		=
		\frac{\gamma\lambda_1}{1+\rho(1-m_0)}.
	\end{equation*}
	Thus \(\alpha=K\) for the published baseline. Define
	\begin{equation*}
		\mathcal L[u]=u_\tau-u_{yy}+\alpha u.
	\end{equation*}
	The \(q\) equation may be written as
	\begin{equation*}
		\mathcal L[q]+R(a,q)-\alpha q=0.
	\end{equation*}
	With embedding parameter \(\varepsilon\in[0,1]\), use the homotopy
	\begin{equation}
		\mathcal L[q]
		+
		\varepsilon\left[R(a,q)-\alpha q\right]
		=0.
		\label{eq:I-homotopy}
	\end{equation}
	At \(\varepsilon=0\) this is the linear reference problem, while
	\(\varepsilon=1\) recovers the original nonlinear \(q\) equation.
	
	Seek expansions
	\begin{equation*}
		q=\sum_{n=0}^{\infty}\varepsilon^n q_n,
		\qquad
		a=\sum_{n=0}^{\infty}\varepsilon^n a_n,
	\end{equation*}
	with the nonhomogeneous data assigned to the zeroth-order term:
	\begin{equation*}
		q_0(0,\tau)=m_0,
		\qquad
		q_n(0,\tau)=0\quad(n\ge1),
	\end{equation*}
	and
	\begin{equation*}
		q_n(y,0)=0,
		\qquad
		q_n(\infty,\tau)=0
		\quad(n\ge0).
	\end{equation*}
	
	Although both fields are written as series for purposes of expanding the
	nonlinear reaction term, Formulation I introduces an HPM hierarchy only
	for \(q\). The coefficients \(a_n\) are subsequently determined from the
	corresponding \(q_n\) through the exact structural reconstruction.
	Expand the nonlinear reaction as
	\begin{equation*}
		R(a,q)=\sum_{n=0}^{\infty}\varepsilon^nR_n.
	\end{equation*}
	Let
	\begin{equation*}
		D_0=1+\rho(1-q_0).
	\end{equation*}
	Since
	\(
	[1+\rho(1-q)]R(a,q)=\lambda_1aq,
	\)
	equating coefficients of equal powers of \(\varepsilon\) gives
	\begin{equation}
		R_n
		=
		\frac{
			\displaystyle
			\lambda_1\sum_{j=0}^{n}a_jq_{n-j}
			+
			\displaystyle
			\rho\sum_{j=1}^{n}q_jR_{n-j}
		}{
			D_0
		},
		\qquad n\ge0,
		\label{eq:I-Rn}
	\end{equation}
	where the second sum vanishes for \(n=0\). In particular,
	\begin{equation*}
		R_0=\frac{\lambda_1a_0q_0}{D_0},
		\qquad
		R_1=
		\frac{\lambda_1(a_0q_1+a_1q_0)+\rho q_1R_0}{D_0}.
	\end{equation*}
	
	Substitution of these series into \eqref{eq:I-homotopy} yields
	\begin{equation*}
		\begin{aligned}
			0={}&
			\mathcal L[q_0]
			+\varepsilon\left(\mathcal L[q_1]+R_0-\alpha q_0\right)\\
			&+\varepsilon^2\left(\mathcal L[q_2]+R_1-\alpha q_1\right)
			+\varepsilon^3\left(\mathcal L[q_3]+R_2-\alpha q_2\right)
			+\cdots .
		\end{aligned}
	\end{equation*}
	Equating coefficients of like powers of \(\varepsilon\) produces the
	linear hierarchy below.
	
	\subsubsection{Zeroth-order problem}
	\label{sec:I-zero}
	
	The coefficient of $\varepsilon^0$ gives
	\begin{equation*}
		\mathcal L[q_0]=0,
	\end{equation*}
	or explicitly,
	\begin{equation*}
		(q_0)_\tau
		=
		(q_0)_{yy}
		-\alpha q_0,
	\end{equation*}
	with
	\begin{equation*}
		q_0(y,0)=0,
		\qquad
		q_0(0,\tau)=m_0,
		\qquad
		q_0(y,\tau)\rightarrow0
		\quad (y\rightarrow\infty).
	\end{equation*}
	
	To solve this zeroth-order problem, take the Laplace transform with respect to
	$\tau$ and define
	\begin{equation*}
		Q_0(y,s)
		=
		\mathcal L_\tau\{q_0(y,\tau)\}.
	\end{equation*}
	Since $q_0(y,0)=0$, the transformed equation is
	\begin{equation*}
		(Q_0)_{yy}
		-
		(s+\alpha)Q_0
		=
		0.
	\end{equation*}
	The solution that remains bounded as $y\rightarrow\infty$ is
	\begin{equation*}
		Q_0(y,s)
		=
		C(s)e^{-y\sqrt{s+\alpha}}.
	\end{equation*}
	Using
	\begin{equation*}
		Q_0(0,s)=\frac{m_0}{s},
	\end{equation*}
	gives
	\begin{equation*}
		Q_0(y,s)
		=
		\frac{m_0}{s}
		e^{-y\sqrt{s+\alpha}}.
	\end{equation*}
	
	Inverting the Laplace transform yields
	\begin{equation}
		q_0(y,\tau)
		=
		\frac{m_0}{2}
		\left[
		e^{-y\sqrt{\alpha}}
		\erfc\left(
		\frac{y}{2\sqrt{\tau}}
		-\sqrt{\alpha\tau}
		\right)
		+
		e^{y\sqrt{\alpha}}
		\erfc\left(
		\frac{y}{2\sqrt{\tau}}
		+\sqrt{\alpha\tau}
		\right)
		\right].
		\label{eq:I-q0closed}
	\end{equation}
	
	Because $\alpha=K$, equation \eqref{eq:I-q0closed} is identical to the
	published catalyst approximation:
	\begin{equation*}
		q_0(y,\tau)=q_{\mathrm P}(y,\tau).
	\end{equation*}
	Thus Formulation I and the published approximation have the same
	zeroth-order catalyst field. Their distinction at this level lies in the
	treatment of the substrate concentration.
	
	Differentiating \eqref{eq:I-q0closed} at the electrode gives
	\begin{equation*}
		-(q_0)_y(0,\tau)
		=
		m_0
		\left[
		\sqrt{\alpha}\,
		\erf\left(\sqrt{\alpha\tau}\right)
		+
		\frac{e^{-\alpha\tau}}{\sqrt{\pi\tau}}
		\right].
	\end{equation*}
	Hence
	\begin{equation*}
		\psi_0(\tau)
		=
		-(q_0)_y(0,\tau)
		=
		\psi_{\mathrm P}(\tau).
	\end{equation*}
	Consequently, improvements in the catalyst concentration and electrode
	current relative to the published baseline arise only after higher-order
	HPM corrections are included.
	
	The substrate field associated with $q_0$ is not taken from the published
	companion approximation. Instead, it is reconstructed from the exact
	structural relation \eqref{eq:arep}:
	\begin{equation*}
		a_0(y,\tau)
		=
		q_0(y,\tau)
		+\gamma
		+
		\int_0^\tau
		\frac{(q_0)_y(0,\sigma)}
		{\sqrt{\pi(\tau-\sigma)}}
		\exp\left[
		-\frac{y^2}{4(\tau-\sigma)}
		\right]
		\,d\sigma.
	\end{equation*}
	Therefore,
	\begin{equation*}
		(a_0-q_0)_\tau
		=
		(a_0-q_0)_{yy},
	\end{equation*}
	exactly.
	
	The zeroth-order reaction coefficient is then
	\begin{equation}
		R_0(y,\tau)
		=
		\frac{
			\lambda_1a_0(y,\tau)q_0(y,\tau)
		}{
			1+\rho\bigl(1-q_0(y,\tau)\bigr)
		}.
		\label{eq:I-R0}
	\end{equation}
	
	Thus the zeroth-order calculation proceeds sequentially as
	\begin{equation*}
		q_0
		\quad\longrightarrow\quad
		a_0
		\quad\longrightarrow\quad
		R_0,
	\end{equation*}
	where $R_0$ provides the nonlinear forcing required for the first-order
	correction.
	
	\subsubsection{First-order problem}
	\label{sec:I-first-order}
	
	Equating coefficients of $\varepsilon^1$ in the homotopy expansion gives
	\begin{equation*}
		\mathcal L[q_1]+R_0-\alpha q_0=0,
	\end{equation*}
	or equivalently,
	\begin{equation*}
		(q_1)_\tau
		=
		(q_1)_{yy}
		-\alpha q_1
		+\alpha q_0-R_0.
	\end{equation*}
	Since $R_0$ is completely determined by $q_0$ and $a_0$ through
	\eqref{eq:I-R0}, define the known first-order forcing
	\begin{equation*}
		F_1(y,\tau)
		=
		\alpha q_0(y,\tau)-R_0(y,\tau).
	\end{equation*}
	The first-order problem is therefore
	\begin{equation*}
		(q_1)_\tau
		=
		(q_1)_{yy}
		-\alpha q_1
		+
		F_1(y,\tau),
		\qquad
		y>0,\quad \tau>0,
	\end{equation*}
	with homogeneous initial and boundary conditions
	\begin{equation*}
		q_1(y,0)=0,
		\qquad
		q_1(0,\tau)=0,
		\qquad
		q_1(y,\tau)\rightarrow0
		\quad (y\rightarrow\infty).
	\end{equation*}
	
	\paragraph{Green-function representation.}
	
	Taking the Laplace transform with respect to $\tau$, let
	\begin{equation*}
		Q_1(y,s)
		=
		\mathcal L_\tau\{q_1(y,\tau)\},
		\qquad
		\widehat F_1(y,s)
		=
		\mathcal L_\tau\{F_1(y,\tau)\}.
	\end{equation*}
	Since $q_1(y,0)=0$, the transformed equation becomes
	\begin{equation*}
		\frac{d^2Q_1}{dy^2}
		-
		(s+\alpha)Q_1
		=
		-\widehat F_1(y,s),
	\end{equation*}
	subject to
	\begin{equation*}
		Q_1(0,s)=0,
		\qquad
		Q_1(y,s)\rightarrow0
		\quad (y\rightarrow\infty).
	\end{equation*}
	
	Setting
	\begin{equation*}
		\mu=\sqrt{s+\alpha},
	\end{equation*}
	the half-line Dirichlet Green function is
	\begin{equation*}
		\widehat G_\alpha(y,\eta,s)
		=
		\frac{1}{2\mu}
		\left[
		e^{-\mu|y-\eta|}
		-
		e^{-\mu(y+\eta)}
		\right],
	\end{equation*}
	where the image term enforces the homogeneous condition at $y=0$.
	Hence
	\begin{equation*}
		Q_1(y,s)
		=
		\int_0^\infty
		\widehat G_\alpha(y,\eta,s)
		\widehat F_1(\eta,s)\,d\eta.
	\end{equation*}
	
	Using
	\begin{equation*}
		\mathcal L^{-1}
		\left\{
		\frac{e^{-b\sqrt{s+\alpha}}}
		{2\sqrt{s+\alpha}}
		\right\}
		=
		\frac{e^{-\alpha t}}{\sqrt{4\pi t}}
		\exp\left(-\frac{b^2}{4t}\right),
	\end{equation*}
	gives the time-domain Green function
	\begin{equation}
		G_\alpha(y,\eta,t)
		=
		\frac{e^{-\alpha t}}{\sqrt{4\pi t}}
		\left[
		\exp\left(-\frac{(y-\eta)^2}{4t}\right)
		-
		\exp\left(-\frac{(y+\eta)^2}{4t}\right)
		\right].
		\label{eq:I-Galpha}
	\end{equation}
	Therefore,
	\begin{equation}
		q_1(y,\tau)
		=
		\int_0^\tau
		\int_0^\infty
		G_\alpha(y,\eta,\tau-\sigma)
		F_1(\eta,\sigma)
		\,d\eta\,d\sigma.
		\label{eq:I-q1green}
	\end{equation}
	Equation \eqref{eq:I-q1green} is an exact semi-analytical representation
	of the first HPM correction.
	
	\paragraph{Regularized time integration.}
	
	The Green kernel in \eqref{eq:I-q1green} contains the endpoint factor
	$(\tau-\sigma)^{-1/2}$. Introduce
	\begin{equation*}
		\sigma=\tau(1-\zeta^2),
		\qquad
		0\leq\zeta\leq1.
	\end{equation*}
	Then
	\begin{equation*}
		\tau-\sigma=\tau\zeta^2,
		\qquad
		d\sigma=-2\tau\zeta\,d\zeta,
	\end{equation*}
	and \eqref{eq:I-q1green} becomes
	\begin{align*}
		q_1(y,\tau)
		={}&
		\frac{\sqrt{\tau}}{\sqrt{\pi}}
		\int_0^1
		e^{-\alpha\tau\zeta^2}
		\int_0^\infty
		\Bigg[
		\exp\left(
		-\frac{(y-\eta)^2}{4\tau\zeta^2}
		\right)
		-
		\exp\left(
		-\frac{(y+\eta)^2}{4\tau\zeta^2}
		\right)
		\Bigg]
		\nonumber\\
		&\hspace{3.0cm}\times
		F_1\left(
		\eta,\tau(1-\zeta^2)
		\right)
		\,d\eta\,d\zeta .
	\end{align*}
	This transformation removes the explicit square-root endpoint factor in
	the time integral and is therefore more convenient for direct quadrature.
	
	\paragraph{Surface derivative and reconstruction of $a_1$.}
	
	Once $q_1$ is known, the first-order contribution to the exact structural
	relation is
	\begin{equation*}
		a_1(y,\tau)
		=
		q_1(y,\tau)
		+
		\int_0^\tau
		\frac{(q_1)_y(0,\sigma)}
		{\sqrt{\pi(\tau-\sigma)}}
		\exp\left[
		-\frac{y^2}{4(\tau-\sigma)}
		\right]
		\,d\sigma.
	\end{equation*}
	
	The required surface derivative can be evaluated directly from the Green
	representation rather than by numerically differentiating $q_1$.
	Differentiating \eqref{eq:I-Galpha} with respect to $y$ and setting
	$y=0$ gives
	\begin{equation*}
		\frac{\partial G_\alpha}{\partial y}(0,\eta,t)
		=
		\frac{\eta e^{-\alpha t}}
		{2\sqrt{\pi}\,t^{3/2}}
		\exp\left(
		-\frac{\eta^2}{4t}
		\right).
	\end{equation*}
	Consequently,
	\begin{equation*}
		(q_1)_y(0,\tau)
		=
		\int_0^\tau
		\int_0^\infty
		\frac{
			\eta e^{-\alpha(\tau-\sigma)}
		}{
			2\sqrt{\pi}(\tau-\sigma)^{3/2}
		}
		\exp\left[
		-\frac{\eta^2}{4(\tau-\sigma)}
		\right]
		F_1(\eta,\sigma)
		\,d\eta\,d\sigma.
	\end{equation*}
	
	Using
	\begin{equation*}
		\sigma=\tau(1-\zeta^2),
		\qquad
		\eta=2\sqrt{\tau}\,\zeta z,
	\end{equation*}
	gives
	\begin{equation*}
		(q_1)_y(0,\tau)
		=
		\frac{4\sqrt{\tau}}{\sqrt{\pi}}
		\int_0^1
		e^{-\alpha\tau\zeta^2}
		\int_0^\infty
		z e^{-z^2}
		F_1\left(
		2\sqrt{\tau}\,\zeta z,
		\tau(1-\zeta^2)
		\right)
		\,dz\,d\zeta.
	\end{equation*}
	This expression provides both the first current correction
	\begin{equation*}
		\psi_1(\tau)=-(q_1)_y(0,\tau)
	\end{equation*}
	and the boundary derivative required in the reconstruction of $a_1$.
	
	Applying the same time transformation to the structural reconstruction of $a_1$ yields
	\begin{equation*}
		a_1(y,\tau)
		=
		q_1(y,\tau)
		+
		\frac{2\sqrt{\tau}}{\sqrt{\pi}}
		\int_0^1
		(q_1)_y\left(
		0,\tau(1-\zeta^2)
		\right)
		\exp\left[
		-\frac{y^2}{4\tau\zeta^2}
		\right]
		\,d\zeta,
	\end{equation*}
	with the endpoint understood through the corresponding limiting value.
	
	The three regularized expressions above provide directly evaluable semi-analytical
	representations of the complete first-order correction.
	
	Once $q_1$ and $a_1$ have been determined, the reaction coefficient
	$R_1$ follows from the recurrence for $R_n$. Thus the first-order
	calculation proceeds as
	\begin{equation*}
		q_0,\ a_0,\ R_0
		\quad\longrightarrow\quad
		q_1
		\quad\longrightarrow\quad
		a_1
		\quad\longrightarrow\quad
		R_1.
	\end{equation*}
	
	No new analytical structure appears at higher orders: the same Green
	kernel generates each catalyst correction, the corresponding substrate
	coefficient is recovered from the exact structural relation, and the next
	reaction coefficient follows from the reaction recurrence. We therefore
	proceed directly to the general $n$th-order formulation.
	
	\subsubsection{General $n$th-order formulation}
	\label{sec:I-general}
	
	The preceding construction extends directly to every higher HPM order.
	For $n\ge1$, equating coefficients of $\varepsilon^n$ gives
	\begin{equation*}
		\mathcal L[q_n]
		+
		R_{n-1}
		-
		\alpha q_{n-1}
		=
		0,
	\end{equation*}
	or equivalently,
	\begin{equation}
		(q_n)_\tau
		=
		(q_n)_{yy}
		-\alpha q_n
		+\alpha q_{n-1}
		-R_{n-1},
		\qquad n\ge1.
		\label{eq:I-qn}
	\end{equation}
	
	All nonhomogeneous initial and boundary data are carried by the
	zeroth-order term. Hence, for every $n\ge1$,
	\begin{equation*}
		q_n(y,0)=0,
		\qquad
		q_n(0,\tau)=0,
		\qquad
		q_n(y,\tau)\rightarrow0
		\quad (y\rightarrow\infty).
	\end{equation*}
	
	Define
	\begin{equation*}
		F_n(y,\tau)
		=
		\alpha q_{n-1}(y,\tau)
		-
		R_{n-1}(y,\tau),
		\qquad n\ge1.
	\end{equation*}
	Then \eqref{eq:I-qn} becomes
	\begin{equation*}
		(q_n)_\tau
		=
		(q_n)_{yy}
		-\alpha q_n
		+
		F_n(y,\tau).
	\end{equation*}
	
	Because the linear operator and homogeneous boundary conditions are the
	same at every order, the Green function $G_\alpha$ in
	\eqref{eq:I-Galpha} applies for all $n\ge1$. Therefore,
	\begin{equation*}
		q_n(y,\tau)
		=
		\int_0^\tau
		\int_0^\infty
		G_\alpha(y,\eta,\tau-\sigma)
		F_n(\eta,\sigma)
		\,d\eta\,d\sigma,
		\qquad n\ge1.
	\end{equation*}
	The regularized representations derived for $q_1$ and its surface
	derivative extend immediately to $q_n$ by replacing $F_1$ with $F_n$.
	
	Once $q_n$ has been obtained, the corresponding substrate coefficient is
	reconstructed from the exact difference relation:
	\begin{equation*}
		a_n(y,\tau)
		=
		q_n(y,\tau)
		+
		\int_0^\tau
		\frac{(q_n)_y(0,\sigma)}
		{\sqrt{\pi(\tau-\sigma)}}
		\exp\left[
		-\frac{y^2}{4(\tau-\sigma)}
		\right]
		\,d\sigma,
		\qquad n\ge1.
	\end{equation*}
	
	The reaction coefficient $R_n$ is then obtained from the recurrence
	\eqref{eq:I-Rn}. Thus the general recursive step may be written
	as
	\begin{equation*}
		\{q_j,a_j,R_j\}_{j=0}^{n-1}
		\quad\longrightarrow\quad
		q_n
		\quad\longrightarrow\quad
		a_n
		\quad\longrightarrow\quad
		R_n,
		\qquad n\ge1.
	\end{equation*}
	Starting from the zeroth-order quantities, the hierarchy has the form
	\begin{equation*}
		q_0
		\longrightarrow
		a_0
		\longrightarrow
		R_0
		\longrightarrow
		q_1
		\longrightarrow
		a_1
		\longrightarrow
		R_1
		\longrightarrow
		q_2
		\longrightarrow
		a_2
		\longrightarrow
		R_2
		\longrightarrow
		\cdots .
	\end{equation*}
	
	After setting the embedding parameter $\varepsilon=1$, the order-$N$
	Formulation-I approximation is
	\begin{equation*}
		q_{\mathrm I}^{[N]}(y,\tau)
		=
		\sum_{n=0}^{N}q_n(y,\tau),
		\qquad
		a_{\mathrm I}^{[N]}(y,\tau)
		=
		\sum_{n=0}^{N}a_n(y,\tau).
	\end{equation*}
	
	Because $a_0$ and every $a_n$, $n\ge1$, are obtained from the
	corresponding catalyst coefficient through the exact heat-equation
	reconstruction, every truncation preserves the structural identity:
	\begin{equation*}
		\left(
		a_{\mathrm I}^{[N]}
		-
		q_{\mathrm I}^{[N]}
		\right)_\tau
		=
		\left(
		a_{\mathrm I}^{[N]}
		-
		q_{\mathrm I}^{[N]}
		\right)_{yy}.
	\end{equation*}
	
	Finally, the corresponding electrode-current approximation is
	\begin{equation*}
		\psi_{\mathrm I}^{[N]}(\tau)
		=
		-
		\sum_{n=0}^{N}
		(q_n)_y(0,\tau).
	\end{equation*}
	
	\subsection{Equivalent Coupled Formulation and Consistency Check}
	\label{sec:coupled-consistency}
	
	The structure-preserving reconstruction in Formulation I can also be
	obtained by applying the same frozen-coefficient homotopy directly to both
	reaction--diffusion equations. This construction is not introduced as a
	separate approximation; rather, it provides an independent derivation of
	the structural reconstruction and a useful consistency check.
	
	Consider
	\begin{equation*}
		q_\tau-q_{yy}+\alpha q
		+
		\varepsilon\bigl[R(a,q)-\alpha q\bigr]
		=0,
	\end{equation*}
	together with
	\begin{equation*}
		a_\tau-a_{yy}+\alpha q
		+
		\varepsilon\bigl[R(a,q)-\alpha q\bigr]
		=0.
	\end{equation*}
	Introduce the expansions
	\begin{equation*}
		q=\sum_{n=0}^{\infty}\varepsilon^n\widetilde q_n,
		\qquad
		a=\sum_{n=0}^{\infty}\varepsilon^n\widetilde a_n,
		\qquad
		R(a,q)=\sum_{n=0}^{\infty}\varepsilon^n\widetilde R_n.
	\end{equation*}
	
	At every order $n\ge1$, the coupled hierarchy has the form
	\begin{equation*}
		(\widetilde q_n)_\tau
		=
		(\widetilde q_n)_{yy}
		-\alpha\widetilde q_n
		+\alpha\widetilde q_{n-1}
		-\widetilde R_{n-1},
	\end{equation*}
	and
	\begin{equation*}
		(\widetilde a_n)_\tau
		=
		(\widetilde a_n)_{yy}
		-\alpha\widetilde q_n
		+\alpha\widetilde q_{n-1}
		-\widetilde R_{n-1}.
	\end{equation*}
	Subtracting gives
	\begin{equation*}
		(\widetilde a_n-\widetilde q_n)_\tau
		=
		(\widetilde a_n-\widetilde q_n)_{yy}.
	\end{equation*}
	
	The difference field above inherits the same
	initial, electrode, and far-field data as the order-$n$ structural
	reconstruction derived from \eqref{eq:arep}. Hence, by uniqueness of the
	corresponding half-line heat problem,
	\begin{equation*}
		\widetilde a_n(y,\tau)
		=
		\widetilde q_n(y,\tau)
		+
		\int_0^\tau
		\frac{(\widetilde q_n)_y(0,\sigma)}
		{\sqrt{\pi(\tau-\sigma)}}
		\exp\left[
		-\frac{y^2}{4(\tau-\sigma)}
		\right]
		\,d\sigma,
		\qquad n\ge1.
	\end{equation*}
	At zeroth order the same argument gives
	\begin{equation*}
		\widetilde a_0(y,\tau)
		=
		\widetilde q_0(y,\tau)
		+\gamma
		+
		\int_0^\tau
		\frac{(\widetilde q_0)_y(0,\sigma)}
		{\sqrt{\pi(\tau-\sigma)}}
		\exp\left[
		-\frac{y^2}{4(\tau-\sigma)}
		\right]
		\,d\sigma.
	\end{equation*}
	
	Thus direct solution of the coupled frozen-coefficient hierarchy produces
	the same substrate coefficients as the structure-preserving reconstruction
	used in Formulation I. Since the catalyst equations and reaction
	recurrence are then identical order by order, the two constructions yield
	the same continuous HPM approximation.
	
	The coupled frozen-coefficient construction therefore does not constitute
	a third formulation. Its role is to provide an independent analytical
	derivation of the reconstruction formula and, later, a numerical
	consistency check on its implementation.
	
	\subsection{Formulation II: Profile-Linearized
		Structure-Preserving HPM}
	\label{sec:formulation-II}
	
	Formulation I uses the constant coefficient $\alpha$ to define its
	auxiliary diffusion--reaction operator. A natural alternative is to retain
	more of the local dependence of the nonlinear reaction term on both
	concentration fields by using a first-order Taylor approximation.
	
	A Taylor expansion about the fixed state
	\(
	(a,q)=(\gamma,m_0)
	\),
	however, is not compatible with the physical far-field state
	\(
	(\gamma,0)
	\).
	Indeed, let
	\begin{equation*}
		D_*
		=
		1+\rho(1-m_0).
	\end{equation*}
	At the fixed reference state,
	\begin{equation*}
		R_*
		=
		\frac{\lambda_1\gamma m_0}{D_*},
		\qquad
		R_{q,*}
		=
		\frac{\lambda_1\gamma(1+\rho)}{D_*^2}.
	\end{equation*}
	Since the $a$-component of the fixed reference state already equals its
	far-field value, the corresponding first-order Taylor approximation
	evaluated at $(\gamma,0)$ is
	\begin{equation*}
		R_{\mathrm{lin},*}(\gamma,0)
		=
		R_*
		-
		m_0R_{q,*}
		=
		-
		\frac{
			\lambda_1\gamma\rho m_0^2
		}{
			D_*^2
		},
	\end{equation*}
	which is generally nonzero. Thus a fixed-state linearization introduces
	an artificial reaction contribution at the physical far field.
	
	To avoid this mismatch, Formulation II linearizes the nonlinear reaction
	term about a space- and time-dependent reference profile that already
	satisfies the initial, electrode, and far-field conditions of the original
	problem.
	
	\subsubsection{Reference profile and local linearization}
	\label{sec:II-reference}
	
	We choose the zeroth-order structure-preserving pair from Formulation I as
	the reference profile:
	\begin{equation*}
		\bar q(y,\tau)=q_0(y,\tau),
		\qquad
		\bar a(y,\tau)=a_0(y,\tau).
	\end{equation*}
	Thus
	\begin{equation*}
		\bar q_\tau
		=
		\bar q_{yy}
		-
		\alpha\bar q,
		\qquad
		\bar a_\tau
		=
		\bar a_{yy}
		-
		\alpha\bar q,
	\end{equation*}
	and consequently
	\begin{equation*}
		(\bar a-\bar q)_\tau
		=
		(\bar a-\bar q)_{yy}.
	\end{equation*}
	
	Because $\bar q=q_0=q_{\mathrm P}$, the reference catalyst profile is the
	same as the published leading approximation. The substrate reference
	profile, however, is the structure-preserving reconstruction $a_0$ rather
	than the published companion approximation $a_{\mathrm P}$.
	
	The reference pair satisfies the original initial conditions,
	\begin{equation*}
		\bar q(y,0)=0,
		\qquad
		\bar a(y,0)=\gamma,
	\end{equation*}
	the electrode conditions,
	\begin{equation*}
		\bar q(0,\tau)=m_0,
		\qquad
		\bar a_y(0,\tau)=0,
	\end{equation*}
	and the physical far-field conditions,
	\begin{equation*}
		\bar q(y,\tau)\rightarrow0,
		\qquad
		\bar a(y,\tau)\rightarrow\gamma,
		\qquad
		y\rightarrow\infty.
	\end{equation*}
	
	Define
	\begin{equation*}
		\bar D
		=
		1+\rho(1-\bar q),
		\qquad
		\bar R
		=
		R(\bar a,\bar q)
		=
		\frac{\lambda_1\bar a\bar q}{\bar D}.
	\end{equation*}
	
	From \eqref{eq:R},
	\begin{equation*}
		R_a(a,q)
		=
		\frac{\lambda_1q}
		{1+\rho(1-q)},
	\end{equation*}
	and
	\begin{equation*}
		R_q(a,q)
		=
		\frac{\lambda_1a(1+\rho)}
		{\left[1+\rho(1-q)\right]^2}.
	\end{equation*}
	Evaluating these derivatives along the reference profile gives
	\begin{equation*}
		\beta_a(y,\tau)
		=
		R_a(\bar a,\bar q)
		=
		\frac{\lambda_1\bar q}{\bar D},
		\qquad
		\beta_q(y,\tau)
		=
		R_q(\bar a,\bar q)
		=
		\frac{\lambda_1\bar a(1+\rho)}{\bar D^2}.
	\end{equation*}
	
	At each point $(y,\tau)$, the first-order Taylor approximation of the
	reaction term about the local reference state
	\(
	(\bar a(y,\tau),\bar q(y,\tau))
	\)
	is therefore
	\begin{equation}
		R_{\mathrm L}(a,q)
		=
		\bar R
		+
		\beta_a(a-\bar a)
		+
		\beta_q(q-\bar q).
		\label{eq:IIp-RL}
	\end{equation}
	
	Unlike a linearization about a fixed concentration state,
	\eqref{eq:IIp-RL} adapts to the evolving reference profile while remaining
	compatible with the physical initial, boundary, and far-field states.
	
	\subsubsection{Exact correction equations}
	\label{sec:II-corrections}
	
	Write the exact solution as the reference profile plus corrections:
	\begin{equation}
		q=\bar q+u,
		\qquad
		a=\bar a+v.
		\label{eq:IIp-corrections}
	\end{equation}
	
	Because the reference profile satisfies
	\begin{equation*}
		\bar q_\tau
		=
		\bar q_{yy}
		-\alpha\bar q,
		\qquad
		\bar a_\tau
		=
		\bar a_{yy}
		-\alpha\bar q,
	\end{equation*}
	its residual with respect to either equation of the original nonlinear
	system is
	\begin{equation*}
		\mathcal E(y,\tau)
		=
		\bar R-\alpha\bar q.
	\end{equation*}
	Indeed,
	\begin{equation*}
		\bar q_\tau-\bar q_{yy}+\bar R
		=
		\bar a_\tau-\bar a_{yy}+\bar R
		=
		\mathcal E.
	\end{equation*}
	Thus $\mathcal E$ measures the error made when the reference pair is
	substituted into the full nonlinear reaction--diffusion system.
	
	Define the nonlinear Taylor remainder by
	\begin{equation*}
		\mathcal N(u,v)
		=
		R(\bar a+v,\bar q+u)
		-\bar R
		-\beta_a v
		-\beta_q u.
	\end{equation*}
	Since
	\begin{equation*}
		1+\rho(1-\bar q-u)
		=
		\bar D-\rho u,
	\end{equation*}
	we have
	\begin{equation*}
		R(\bar a+v,\bar q+u)
		=
		\frac{
			\lambda_1(\bar a+v)(\bar q+u)
		}{
			\bar D-\rho u
		}.
	\end{equation*}
	Substituting the definitions of $\bar R$, $\beta_a$, and $\beta_q$ into this remainder and simplifying gives the exact expression
	\begin{equation}
		\mathcal N(u,v)
		=
		\frac{
			\lambda_1(1+\rho)\,
			u\left(
			\bar D\,v+\rho\bar a\,u
			\right)
		}{
			\bar D^2(\bar D-\rho u)
		}.
		\label{eq:IIp-Nexact}
	\end{equation}
	
	Provided $\bar D-\rho u$ remains bounded away from zero,
	\eqref{eq:IIp-Nexact} shows that
	\begin{equation*}
		\mathcal N(u,v)
		=
		O\!\left(u^2+|uv|\right)
		\qquad
		\text{as }(u,v)\to(0,0).
	\end{equation*}
	Thus all terms that are linear in the corrections are contained in the
	profile-dependent Taylor operator. No $v^2$ term occurs because the
	reaction function $R(a,q)$ is linear in $a$ for fixed $q$.
	
	Substituting \eqref{eq:IIp-corrections} into the original nonlinear
	system and using
	\begin{equation*}
		R(\bar a+v,\bar q+u)
		=
		\bar R+\beta_a v+\beta_q u+\mathcal N(u,v)
	\end{equation*}
	gives the exact correction system
	\begin{align}
		u_\tau
		&=
		u_{yy}
		-\beta_q u
		-\beta_a v
		-\mathcal E
		-\mathcal N(u,v),
		\nonumber\\
		v_\tau
		&=
		v_{yy}
		-\beta_q u
		-\beta_a v
		-\mathcal E
		-\mathcal N(u,v).
		\label{eq:IIp-correction-system}
	\end{align}
	
	Because the reference profile satisfies all nonhomogeneous data of the
	original problem, the corrections satisfy homogeneous initial, electrode,
	and far-field conditions:
	\begin{equation*}
		u(y,0)=0,
		\qquad
		u(0,\tau)=0,
		\qquad
		u(y,\tau)\rightarrow0
		\quad (y\rightarrow\infty),
	\end{equation*}
	and
	\begin{equation*}
		v(y,0)=0,
		\qquad
		v_y(0,\tau)=0,
		\qquad
		v(y,\tau)\rightarrow0
		\quad (y\rightarrow\infty).
	\end{equation*}
	
	\subsubsection{Homotopy hierarchy}
	\label{sec:II-homotopy}
	
	The exact correction system \eqref{eq:IIp-correction-system} contains two
	distinct contributions beyond the profile-dependent linear operator. The
	known function $\mathcal E$ is the residual of the reference profile,
	whereas $\mathcal N(u,v)$ contains the nonlinear terms remaining after
	the first-order Taylor expansion. We therefore introduce the embedding
	parameter only in the nonlinear remainder:
	\begin{equation}
		\begin{aligned}
			u_\tau-u_{yy}
			+\beta_q u+\beta_a v+\mathcal E
			+\varepsilon\,\mathcal N(u,v)
			&=0,
			\\
			v_\tau-v_{yy}
			+\beta_q u+\beta_a v+\mathcal E
			+\varepsilon\,\mathcal N(u,v)
			&=0,
		\end{aligned}
		\qquad
		0\leq\varepsilon\leq1.
		\label{eq:IIp-homotopy}
	\end{equation}
	
	At $\varepsilon=0$, the nonlinear Taylor remainder is suppressed, but
	the residual of the reference profile is retained. Thus the leading
	correction solves the profile-linearized response to the known forcing
	$-\mathcal E$. At $\varepsilon=1$, equation
	\eqref{eq:IIp-homotopy} recovers the exact nonlinear correction system
	\eqref{eq:IIp-correction-system}.
	
	Expand the correction fields as
	\begin{equation*}
		u
		=
		\sum_{n=0}^{\infty}\varepsilon^n u_n,
		\qquad
		v
		=
		\sum_{n=0}^{\infty}\varepsilon^n v_n,
	\end{equation*}
	and, after substituting these series into the nonlinear remainder, write
	\begin{equation}
		\mathcal N(u,v)
		=
		\sum_{n=0}^{\infty}
		\varepsilon^n\mathcal N_n.
		\label{eq:IIp-Nseries}
	\end{equation}
	Here $\mathcal N_n$ denotes the coefficient of $\varepsilon^n$ generated
	by the expansions of $u$ and $v$; an explicit recurrence for these
	coefficients is derived below.
	
	At order $\varepsilon^0$, the leading correction pair
	$(u_0,v_0)$ satisfies
	\begin{equation}
		\begin{aligned}
			(u_0)_\tau
			&=
			(u_0)_{yy}
			-\beta_q u_0
			-\beta_a v_0
			-\mathcal E,
			\\
			(v_0)_\tau
			&=
			(v_0)_{yy}
			-\beta_q u_0
			-\beta_a v_0
			-\mathcal E.
		\end{aligned}
		\label{eq:IIp-order0}
	\end{equation}
	Thus $(u_0,v_0)$ is the profile-linearized correction generated by the
	residual of the reference pair.
	
	Because
	\begin{equation*}
		\varepsilon\mathcal N(u,v)
		=
		\varepsilon\mathcal N_0
		+
		\varepsilon^2\mathcal N_1
		+
		\varepsilon^3\mathcal N_2
		+\cdots,
	\end{equation*}
	the coefficient of $\varepsilon^n$, for $n\ge1$, gives
	\begin{equation*}
		\begin{aligned}
			(u_n)_\tau
			&=
			(u_n)_{yy}
			-\beta_q u_n
			-\beta_a v_n
			-\mathcal N_{n-1},
			\\
			(v_n)_\tau
			&=
			(v_n)_{yy}
			-\beta_q u_n
			-\beta_a v_n
			-\mathcal N_{n-1}.
		\end{aligned}
	\end{equation*}
	
	Every correction pair satisfies homogeneous data:
	\begin{equation*}
		u_n(y,0)=0,
		\qquad
		v_n(y,0)=0,
	\end{equation*}
	\begin{equation*}
		u_n(0,\tau)=0,
		\qquad
		(v_n)_y(0,\tau)=0,
	\end{equation*}
	and
	\begin{equation*}
		u_n(y,\tau)\rightarrow0,
		\qquad
		v_n(y,\tau)\rightarrow0,
		\qquad
		y\rightarrow\infty.
	\end{equation*}
	
	Consequently, the same coupled variable-coefficient linear operator
	appears at every HPM level. The leading correction is forced by
	$-\mathcal E$, while each subsequent correction is forced by
	$-\mathcal N_{n-1}$. The nonlinear calculation therefore reduces to a
	sequence of linear correction problems with a common operator and
	successively determined forcing terms.
	
	\subsubsection{Leading profile correction and integral representation}
	\label{sec:II-leading-correction}
	
	With no profile correction retained, Formulation II reduces to the
	reference pair,
	\begin{equation*}
		q_{\mathrm{II}}^{[0]}=\bar q=q_0,
		\qquad
		a_{\mathrm{II}}^{[0]}=\bar a=a_0.
	\end{equation*}
	Consequently,
	\begin{equation*}
		\psi_{\mathrm{II}}^{[0]}(\tau)
		=
		-\bar q_y(0,\tau)
		=
		\psi_0(\tau)
		=
		\psi_{\mathrm P}(\tau).
	\end{equation*}
	Thus the first distinction between the catalyst and current predictions of
	Formulations I and II appears only after correction terms are included.
	
	The leading correction pair $(u_0,v_0)$ satisfies
	\eqref{eq:IIp-order0}. Because the coefficients
	$\beta_a(y,\tau)$ and $\beta_q(y,\tau)$ depend on the evolving reference
	profile, this is a coupled variable-coefficient system and does not in
	general admit the constant-coefficient error-function representation
	available for $\bar q$.
	
	An exact integral-equation representation can nevertheless be obtained
	using the half-line heat kernels
	\begin{equation*}
		G_D(y,\eta,t)
		=
		\frac{1}{\sqrt{4\pi t}}
		\left[
		\exp\left(
		-\frac{(y-\eta)^2}{4t}
		\right)
		-
		\exp\left(
		-\frac{(y+\eta)^2}{4t}
		\right)
		\right],
	\end{equation*}
	and
	\begin{equation*}
		G_N(y,\eta,t)
		=
		\frac{1}{\sqrt{4\pi t}}
		\left[
		\exp\left(
		-\frac{(y-\eta)^2}{4t}
		\right)
		+
		\exp\left(
		-\frac{(y+\eta)^2}{4t}
		\right)
		\right],
	\end{equation*}
	which enforce homogeneous Dirichlet and Neumann conditions,
	respectively, at $y=0$.
	
	Define
	\begin{equation*}
		\Phi_0(y,\tau)
		=
		\beta_q(y,\tau)u_0(y,\tau)
		+
		\beta_a(y,\tau)v_0(y,\tau)
		+
		\mathcal E(y,\tau).
	\end{equation*}
	Then \eqref{eq:IIp-order0} can be written as
	\begin{equation*}
		(u_0)_\tau-(u_0)_{yy}
		=
		-\Phi_0,
		\qquad
		(v_0)_\tau-(v_0)_{yy}
		=
		-\Phi_0,
	\end{equation*}
	and Duhamel's principle gives
	\begin{equation*}
		u_0(y,\tau)
		=
		-
		\int_0^\tau
		\int_0^\infty
		G_D(y,\eta,\tau-\sigma)
		\Phi_0(\eta,\sigma)
		\,d\eta\,d\sigma,
	\end{equation*}
	and
	\begin{equation*}
		v_0(y,\tau)
		=
		-
		\int_0^\tau
		\int_0^\infty
		G_N(y,\eta,\tau-\sigma)
		\Phi_0(\eta,\sigma)
		\,d\eta\,d\sigma.
	\end{equation*}
	
	Unlike the Green-function representation for $q_1$ in Formulation I,
	these two Green-function representations are not
	direct quadrature formulas involving only previously known functions.
	The source $\Phi_0$ contains the unknown correction fields
	$u_0$ and $v_0$ themselves. They therefore form a coupled Volterra
	integral system equivalent to the variable-coefficient linear problem
	\eqref{eq:IIp-order0}.
	
	The exact difference structure provides an alternative representation of
	$v_0$. Subtracting the two equations in \eqref{eq:IIp-order0} gives
	\begin{equation*}
		(v_0-u_0)_\tau
		=
		(v_0-u_0)_{yy}.
	\end{equation*}
	Moreover,
	\begin{equation*}
		(v_0-u_0)(y,0)=0,
		\qquad
		(v_0-u_0)(y,\tau)\rightarrow0
		\quad (y\rightarrow\infty),
	\end{equation*}
	and, at the electrode,
	\begin{equation*}
		(v_0-u_0)_y(0,\tau)
		=
		-(u_0)_y(0,\tau).
	\end{equation*}
	Hence
	\begin{equation*}
		v_0(y,\tau)
		=
		u_0(y,\tau)
		+
		\int_0^\tau
		\frac{(u_0)_y(0,\sigma)}
		{\sqrt{\pi(\tau-\sigma)}}
		\exp\left[
		-\frac{y^2}{4(\tau-\sigma)}
		\right]
		\,d\sigma.
	\end{equation*}
	
	The first corrected approximation is therefore
	\begin{equation*}
		q_{\mathrm{II}}^{[1]}
		=
		\bar q+u_0,
		\qquad
		a_{\mathrm{II}}^{[1]}
		=
		\bar a+v_0,
	\end{equation*}
	where $v_0$ may be obtained either from the coupled correction problem or
	from the equivalent structural reconstruction above. The corresponding correction to the
	electrode current is $-(u_0)_y(0,\tau)$.
	
	\subsubsection{Nonlinear coefficient recursion}
	\label{sec:II-nonlinear-recursion}
	
	To advance the HPM hierarchy beyond the leading correction, the
	coefficients $\mathcal N_n$ in
	\eqref{eq:IIp-Nseries} must be determined from the previously computed
	correction fields. Rather than re-expanding the full reaction function,
	we use the exact remainder formula \eqref{eq:IIp-Nexact} directly.
	
	Write
	\begin{equation*}
		u
		=
		\sum_{n=0}^{\infty}\varepsilon^n u_n,
		\qquad
		v
		=
		\sum_{n=0}^{\infty}\varepsilon^n v_n,
		\qquad
		\mathcal N
		=
		\sum_{n=0}^{\infty}\varepsilon^n\mathcal N_n.
	\end{equation*}
	Define
	\begin{equation*}
		P(u,v)
		=
		\frac{\lambda_1(1+\rho)}{\bar D^2}
		u\left(
		\bar D\,v+\rho\bar a\,u
		\right),
	\end{equation*}
	so that
	\begin{equation*}
		\mathcal N(u,v)
		=
		\frac{P(u,v)}
		{\bar D-\rho u}.
	\end{equation*}
	
	Expanding
	\begin{equation*}
		P(u,v)
		=
		\sum_{n=0}^{\infty}\varepsilon^n P_n,
	\end{equation*}
	the coefficient of $\varepsilon^n$ is
	\begin{equation*}
		P_n
		=
		\frac{\lambda_1(1+\rho)}{\bar D^2}
		\left[
		\bar D
		\sum_{j=0}^{n}u_jv_{n-j}
		+
		\rho\bar a
		\sum_{j=0}^{n}u_ju_{n-j}
		\right].
	\end{equation*}
	
	Since
	\begin{equation*}
		(\bar D-\rho u)\mathcal N=P,
	\end{equation*}
	equating coefficients of equal powers of $\varepsilon$ gives, at
	zeroth order,
	\begin{equation*}
		\mathcal N_0
		=
		\frac{P_0}
		{\bar D-\rho u_0}
		=
		\frac{
			\lambda_1(1+\rho)
			u_0\left(
			\bar D\,v_0+\rho\bar a\,u_0
			\right)
		}{
			\bar D^2(\bar D-\rho u_0)
		}.
	\end{equation*}
	For every $n\ge1$,
	\begin{equation*}
		\mathcal N_n
		=
		\frac{
			P_n
			+
			\rho\displaystyle\sum_{j=1}^{n}
			u_j\mathcal N_{n-j}
		}{
			\bar D-\rho u_0
		}.
	\end{equation*}
	
	Thus each nonlinear coefficient is computed only after the corresponding
	correction pair has been obtained. In particular, the recursive sequence
	is
	\begin{equation*}
		(u_0,v_0)
		\longrightarrow
		\mathcal N_0
		\longrightarrow
		(u_1,v_1)
		\longrightarrow
		\mathcal N_1
		\longrightarrow
		(u_2,v_2)
		\longrightarrow
		\cdots .
	\end{equation*}
	Hence the forcing $-\mathcal N_{n-1}$ in the $n$th correction problem is
	known before $(u_n,v_n)$ is solved, and no nonlinear solve is required
	within an individual HPM level.

	\subsubsection{Structure preservation and far-field compatibility}
	\label{sec:II-structure}
	
	Because the two homotopy equations in \eqref{eq:IIp-homotopy} contain
	identical reaction, residual, and nonlinear-remainder terms, subtracting
	them gives, for every $\varepsilon$,
	\begin{equation*}
		(v-u)_\tau
		=
		(v-u)_{yy}.
	\end{equation*}
	Using
	\begin{equation*}
		v-u
		=
		\sum_{n=0}^{\infty}
		\varepsilon^n(v_n-u_n),
	\end{equation*}
	and equating coefficients of equal powers of $\varepsilon$ yields
	\begin{equation}
		(v_n-u_n)_\tau
		=
		(v_n-u_n)_{yy},
		\qquad n\ge0.
		\label{eq:IIp-structure}
	\end{equation}
	
	For each $n\ge0$, let
	\begin{equation*}
		d_n=v_n-u_n.
	\end{equation*}
	The homogeneous correction data imply
	\begin{equation*}
		d_n(y,0)=0,
		\qquad
		d_n(y,\tau)\rightarrow0
		\quad (y\rightarrow\infty),
	\end{equation*}
	while at the electrode,
	\begin{equation*}
		(d_n)_y(0,\tau)
		=
		(v_n)_y(0,\tau)
		-
		(u_n)_y(0,\tau)
		=
		-(u_n)_y(0,\tau).
	\end{equation*}
	Thus the same half-line heat-equation representation used in
	Section~\ref{sec:structural} gives
	\begin{equation*}
		v_n(y,\tau)
		=
		u_n(y,\tau)
		+
		\int_0^\tau
		\frac{(u_n)_y(0,\sigma)}
		{\sqrt{\pi(\tau-\sigma)}}
		\exp\left[
		-\frac{y^2}{4(\tau-\sigma)}
		\right]
		\,d\sigma,
		\qquad n\ge0.
	\end{equation*}
	Consequently, each substrate correction may either be solved directly
	from the coupled linear problem or reconstructed from the corresponding
	catalyst correction.
	
	The profile-dependent linearization is also compatible with the physical
	far-field state. As $y\rightarrow\infty$,
	\begin{equation*}
		\bar q\rightarrow0,
		\qquad
		\bar a\rightarrow\gamma,
		\qquad
		\bar D\rightarrow1+\rho,
		\qquad
		\bar R\rightarrow0.
	\end{equation*}
	Since
	\begin{equation*}
		\mathcal E
		=
		\bar R-\alpha\bar q,
	\end{equation*}
	it follows that
	\begin{equation*}
		\mathcal E\rightarrow0.
	\end{equation*}
	Moreover,
	\begin{equation*}
		\beta_a
		=
		\frac{\lambda_1\bar q}{\bar D}
		\rightarrow0,
		\qquad
		\beta_q
		=
		\frac{\lambda_1\bar a(1+\rho)}{\bar D^2}
		\rightarrow
		\frac{\lambda_1\gamma}{1+\rho}.
	\end{equation*}
	
	Finally, at the physical far-field state
	\(
	(a,q)=(\gamma,0)
	\),
	the reference state approaches the same limit,
	\(
	(\bar a,\bar q)\rightarrow(\gamma,0)
	\).
	Therefore,
	\begin{align*}
		R_{\mathrm L}(\gamma,0)
		&=
		\bar R
		+
		\beta_a(\gamma-\bar a)
		+
		\beta_q(0-\bar q)
		\\
		&\longrightarrow0.
	\end{align*}
	Hence the profile-linearized operator introduces no artificial
	far-field reaction forcing, in contrast to a Taylor linearization about
	the fixed state $(\gamma,m_0)$.
	
	\subsubsection{Truncated approximation}
	\label{sec:II-truncation}
	
	We define the truncation index $N$ to denote the number of
	profile-correction levels retained beyond the reference pair. Thus the
	uncorrected reference approximation is
	\begin{equation*}
		q_{\mathrm{II}}^{[0]}(y,\tau)
		=
		\bar q(y,\tau),
		\qquad
		a_{\mathrm{II}}^{[0]}(y,\tau)
		=
		\bar a(y,\tau).
	\end{equation*}
	
	For $N\ge1$, retaining the correction pairs
	\[
	(u_0,v_0),\ldots,(u_{N-1},v_{N-1})
	\]
	gives
	\begin{equation*}
		q_{\mathrm{II}}^{[N]}(y,\tau)
		=
		\bar q(y,\tau)
		+
		\sum_{n=0}^{N-1}u_n(y,\tau),
		\qquad
		a_{\mathrm{II}}^{[N]}(y,\tau)
		=
		\bar a(y,\tau)
		+
		\sum_{n=0}^{N-1}v_n(y,\tau).
	\end{equation*}
	
	In particular,
	\begin{equation*}
		q_{\mathrm{II}}^{[1]}
		=
		\bar q+u_0,
		\qquad
		a_{\mathrm{II}}^{[1]}
		=
		\bar a+v_0,
	\end{equation*}
	so the superscript $[N]$ counts retained profile-correction levels rather
	than the highest correction index.
	
	Since the dimensionless current is determined by the surface gradient of
	$q$, the corresponding current approximation is
	\begin{equation*}
		\psi_{\mathrm{II}}^{[N]}(\tau)
		=
		-\bar q_y(0,\tau)
		-
		\sum_{n=0}^{N-1}
		(u_n)_y(0,\tau),
		\qquad N\ge1,
	\end{equation*}
	with
	\begin{equation*}
		\psi_{\mathrm{II}}^{[0]}(\tau)
		=
		-\bar q_y(0,\tau).
	\end{equation*}
	
	By \eqref{eq:IIp-structure}, every truncated pair
	\(
	(q_{\mathrm{II}}^{[N]},a_{\mathrm{II}}^{[N]})
	\)
	preserves the exact difference heat equation.
	
	\subsection{Summary of the Published and Proposed Solutions}
	\label{sec:hpm-summary}
	
	The principal solution representations for the published baseline and the
	two proposed structure-preserving HPM formulations are summarized in
	Table~\ref{tab:general-solutions}. Detailed derivations of the reference
	solutions, Green functions, reaction coefficients, structural
	reconstructions, and correction equations are given in the preceding
	subsections.
	
	\renewcommand{\arraystretch}{1.35}
	
	\begin{longtable}{
			p{0.17\textwidth}
			p{0.10\textwidth}
			p{0.66\textwidth}}
		
		\caption{
			General solution representations for the published baseline and the two
			structure-preserving HPM formulations.
		}
		\label{tab:general-solutions}
		\\
		
		\toprule
		\textbf{Method}
		&
		\textbf{Quantity}
		&
		\textbf{General solution}
		\\
		\midrule
		\endfirsthead
		
		\multicolumn{3}{c}
		{\tablename\ \thetable\ -- continued from previous page}
		\\
		\toprule
		\textbf{Method}
		&
		\textbf{Quantity}
		&
		\textbf{General solution}
		\\
		\midrule
		\endhead
		
		\midrule
		\multicolumn{3}{r}{Continued on next page}
		\\
		\endfoot
		
		\bottomrule
		\endlastfoot
		
		
		\textbf{Published baseline}
		&
		$q$
		&
		\(
		q_{\mathrm P}(y,\tau)
		=
		\frac{m_0}{2}
		\left[
		e^{-y\sqrt K}
		\erfc\left(
		\frac{y}{2\sqrt{\tau}}-\sqrt{K\tau}
		\right)
		+
		e^{y\sqrt K}
		\erfc\left(
		\frac{y}{2\sqrt{\tau}}+\sqrt{K\tau}
		\right)
		\right].
		\)
		\\[1ex]
		
		&
		$a$
		&
		\(
		a_{\mathrm P}(y,\tau)
		=
		q_{\mathrm P}(y,\tau)+\gamma
		-
		m_0
		\erfc\left(
		\frac{y}{2\sqrt{\tau}}
		\right)
		\left[
		\sqrt{\pi K\tau}\,
		\erf(\sqrt{K\tau})
		+
		e^{-K\tau}
		\right].
		\)
		\\[1ex]
		
		&
		$\psi$
		&
		\(
		\psi_{\mathrm P}(\tau)
		=
		-q_{\mathrm P,y}(0,\tau)
		=
		m_0
		\left[
		\sqrt K\,
		\erf(\sqrt{K\tau})
		+
		\frac{e^{-K\tau}}{\sqrt{\pi\tau}}
		\right].
		\)
		\\
		
		\midrule
		
		
		\textbf{Formulation I}
		&
		$q$
		&
		\(
		q_{\mathrm I}^{[N]}(y,\tau)
		=
		\displaystyle\sum_{n=0}^{N}q_n(y,\tau).
		\)
		\\
		&&
		\(
		q_0(y,\tau)
		=
		\frac{m_0}{2}
		\left[
		e^{-y\sqrt{\alpha}}
		\erfc\left(
		\frac{y}{2\sqrt{\tau}}-\sqrt{\alpha\tau}
		\right)
		+
		e^{y\sqrt{\alpha}}
		\erfc\left(
		\frac{y}{2\sqrt{\tau}}+\sqrt{\alpha\tau}
		\right)
		\right].
		\)
		\\
		&&
		For $n\ge1$,
		\(
		q_n(y,\tau)
		=
		\displaystyle
		\int_0^\tau
		\int_0^\infty
		G_\alpha(y,\eta,\tau-\sigma)
		F_n(\eta,\sigma)
		\,d\eta\,d\sigma,
		\)
		where
		\\
		&&
		\(
		F_n
		=
		\alpha q_{n-1}-R_{n-1},
		\)
		and
		\(
		G_\alpha(y,\eta,t)
		=
		\frac{e^{-\alpha t}}{\sqrt{4\pi t}}
		\left[
		e^{-(y-\eta)^2/(4t)}
		-
		e^{-(y+\eta)^2/(4t)}
		\right].
		\)
		\\[1ex]
		
		&
		$a$
		&
		\(
		a_{\mathrm I}^{[N]}(y,\tau)
		=
		\displaystyle\sum_{n=0}^{N}a_n(y,\tau).
		\)
		\\
		&&
		\(
		a_0(y,\tau)
		=
		q_0(y,\tau)
		+\gamma
		+\mathcal H[q_0](y,\tau).
		\)
		\\
		&&
		For $n\ge1$,
		\(
		a_n(y,\tau)
		=
		q_n(y,\tau)
		+\mathcal H[q_n](y,\tau),
		\)
		where
		\\
		&&
		\(
		\mathcal H[f](y,\tau)
		=
		\displaystyle
		\int_0^\tau
		\frac{f_y(0,\sigma)}
		{\sqrt{\pi(\tau-\sigma)}}
		\exp\left[
		-\frac{y^2}{4(\tau-\sigma)}
		\right]
		\,d\sigma.
		\)
		\\[1ex]
		
		&
		$R_n$
		&
		\(
		R_n
		=
		\frac{
			\displaystyle
			\lambda_1\sum_{j=0}^{n}a_jq_{n-j}
			+
			\displaystyle
			\rho\sum_{j=1}^{n}q_jR_{n-j}
		}{
			1+\rho(1-q_0)
		},
		n\ge0,
		\)
		where the second sum vanishes for $n=0$.
		\\[1ex]
		
		&
		$\psi$
		&
		\(
		\psi_{\mathrm I}^{[N]}(\tau)
		=
		-
		\displaystyle\sum_{n=0}^{N}
		(q_n)_y(0,\tau).
		\)
		\\
		
		\midrule
		
		
		\textbf{Formulation II}
		&
		Reference
		&
		The reference profile is the zeroth-order structure-preserving pair from
		Formulation I:
		\(
		\bar q=q_0=q_{\mathrm P},
		\qquad
		\bar a=a_0.
		\)
		Define
		\\
		&&
		\(
		\bar D
		=
		1+\rho(1-\bar q),
		\qquad
		\bar R
		=
		\frac{\lambda_1\bar a\bar q}{\bar D},
		\)
		\\
		&&
		\(
		\beta_a
		=
		\frac{\lambda_1\bar q}{\bar D},
		\qquad
		\beta_q
		=
		\frac{\lambda_1\bar a(1+\rho)}{\bar D^2},
		\qquad
		\mathcal E
		=
		\bar R-\alpha\bar q.
		\)
		\\[1ex]
		
		&
		$q$
		&
		\(
		q_{\mathrm{II}}^{[0]}(y,\tau)
		=
		\bar q(y,\tau).
		\)
		\\
		&&
		For $N\ge1$,
		\(
		q_{\mathrm{II}}^{[N]}(y,\tau)
		=
		\bar q(y,\tau)
		+
		\displaystyle\sum_{n=0}^{N-1}u_n(y,\tau).
		\)
		\\
		&&
		The leading correction satisfies
		\(
		(u_0)_\tau
		=
		(u_0)_{yy}
		-\beta_q u_0
		-\beta_a v_0
		-\mathcal E.
		\)
		\\
		&&
		For $n\ge1$,
		\(
		(u_n)_\tau
		=
		(u_n)_{yy}
		-\beta_q u_n
		-\beta_a v_n
		-\mathcal N_{n-1}.
		\)
		\\[1ex]
		
		&
		$a$
		&
		\(
		a_{\mathrm{II}}^{[0]}(y,\tau)
		=
		\bar a(y,\tau).
		\)
		\\
		&&
		For $N\ge1$,
		\(
		a_{\mathrm{II}}^{[N]}(y,\tau)
		=
		\bar a(y,\tau)
		+
		\displaystyle\sum_{n=0}^{N-1}v_n(y,\tau).
		\)
		\\
		&&
		The leading correction satisfies
		\(
		(v_0)_\tau
		=
		(v_0)_{yy}
		-\beta_q u_0
		-\beta_a v_0
		-\mathcal E.
		\)
		\\
		&&
		For $n\ge1$,
		\(
		(v_n)_\tau
		=
		(v_n)_{yy}
		-\beta_q u_n
		-\beta_a v_n
		-\mathcal N_{n-1}.
		\)
		\\
		&&
		Equivalently, for every $n\ge0$,
		\(
		v_n(y,\tau)
		=
		u_n(y,\tau)
		+
		\mathcal H[u_n](y,\tau).
		\)
		\\[1ex]
		
		&
		$\mathcal N_n$
		&
		The nonlinear Taylor remainder is
		\(
		\mathcal N(u,v)
		=
		R(\bar a+v,\bar q+u)
		-\bar R-\beta_a v-\beta_q u,
		\)
		or equivalently
		\\
		&&
		\(
		\mathcal N(u,v)
		=
		\frac{
			\lambda_1(1+\rho)\,
			u\left(
			\bar D\,v+\rho\bar a\,u
			\right)
		}{
			\bar D^2(\bar D-\rho u)
		}.
		\)
		\\
		&&
		Writing
		\(
		\mathcal N
		=
		\displaystyle\sum_{n=0}^{\infty}
		\varepsilon^n\mathcal N_n
		\)
		and defining
		\[
		P_n
		=
		\frac{\lambda_1(1+\rho)}{\bar D^2}
		\left[
		\bar D
		\sum_{j=0}^{n}u_jv_{n-j}
		+
		\rho\bar a
		\sum_{j=0}^{n}u_ju_{n-j}
		\right],
		\]
		gives
		\[
		\mathcal N_0
		=
		\frac{P_0}{\bar D-\rho u_0},
		\]
		and, for $n\ge1$,
		\[
		\mathcal N_n
		=
		\frac{
			P_n
			+
			\rho\displaystyle\sum_{j=1}^{n}
			u_j\mathcal N_{n-j}
		}{
			\bar D-\rho u_0
		}.
		\]
		\\[1ex]
		
		&
		$\psi$
		&
		At the uncorrected reference level,
		\(
		\psi_{\mathrm{II}}^{[0]}(\tau)
		=
		-\bar q_y(0,\tau).
		\)
		\\
		&&
		For $N\ge1$,
		\(
		\psi_{\mathrm{II}}^{[N]}(\tau)
		=
		-\bar q_y(0,\tau)
		-
		\displaystyle\sum_{n=0}^{N-1}
		(u_n)_y(0,\tau).
		\)
		\\
		
	\end{longtable}
	
	The published baseline and both proposed formulations share the same
	leading catalyst profile,
	\begin{equation*}
		q_{\mathrm P}=q_0=\bar q,
	\end{equation*}
	because the published coefficient $K$ coincides with the reference
	coefficient $\alpha$. Their differences arise from the treatment of the
	substrate field and from the construction of higher-order corrections.
	
	Formulation I retains a scalar constant-coefficient Green-function
	hierarchy for the catalyst field and reconstructs each substrate
	coefficient from the exact difference relation. Formulation II instead
	uses a profile-dependent Taylor operator and generates a sequence of
	coupled variable-coefficient corrections while preserving the same exact
	difference structure.
	
	The coupled frozen-coefficient construction discussed in
	Section~\ref{sec:coupled-consistency} is not a third approximation. When
	the same auxiliary coefficient $\alpha$ is used, it is mathematically
	equivalent to Formulation I and is retained only as an independent
	analytical and numerical consistency check.
	
	The analytical formulations are now complete. The next section describes
	their numerical evaluation, the independent nonlinear PDE reference
	solution, the computation of the electrode current, and the error and
	residual measures used in the comparisons.

	\section{Numerical Methodology}
	\label{sec:numerical-methodology}
	
	Numerical simulation has long played an important role in
	electrochemical reaction--diffusion analysis, including cyclic-voltammetric
	reaction models, stability analysis of finite-difference schemes, high-order
	boundary-gradient evaluation, and high-order spatial discretization
	\cite{VynnyckyMcKeeBieniasz2021}.
	The analytical formulations in Section~\ref{sec:hpm} lead here to two
	related but numerically distinct tasks. Formulation~I has a scalar
	constant-coefficient correction hierarchy for the catalyst concentration,
	whereas Formulation~II requires the solution of coupled
	variable-coefficient correction systems. In both cases, the exact structural
	relation between the catalyst and substrate fields is retained. An
	independent numerical solution of the original nonlinear reaction--diffusion
	system is also computed as a reference benchmark.
	
	For the production calculations reported in
	Section~\ref{sec:numerical-results}, all grid-based linear and nonlinear
	problems are discretized using fourth-order finite differences in space.
	High-order spatial discretizations are particularly useful in
	electrochemical reaction--diffusion calculations because both concentration
	profiles and boundary gradients must be resolved accurately
	\cite{Bieniasz2003,Bieniasz2004}.
	The resulting semidiscrete systems are integrated in time using the adaptive
	stiff solver \texttt{ode15s}. Analytical Green-function and heat-kernel
	representations are retained where they provide an independent validation
	or a more accurate evaluation of the electrode flux.
	
	The semi-infinite spatial interval is truncated to
	\[
	0\leq y\leq L.
	\]
	A uniform spatial grid is defined by
	\(
	y_j=j\Delta y,
	j=0,\ldots,J,
	\Delta y=\frac{L}{J},
	\)
	with
	\(
	N_y=J+1.
	\)
	The solutions are stored at the uniformly spaced output times
	\[
	\tau^k=k\Delta\tau,
	\qquad
	k=0,\ldots,K,
	\qquad
	\Delta\tau=\frac{\tau_{\max}}{K},
	\]
	so that
	\(
	N_t=K+1.
	\)
	The use of a uniform output-time grid does not imply fixed-step time
	integration: \texttt{ode15s} selects its internal time steps adaptively and
	returns the solution at the prescribed output times.
	
	\subsection{Fourth-Order Spatial Discretization}
	\label{sec:num-spatial}
	
	High-order finite-difference discretizations have been developed
	specifically for electrochemical reaction--diffusion problems
	\cite{Bieniasz2004}. For interior grid points sufficiently far from the
	boundaries, the second derivative is approximated by the fourth-order
	centered formula
	\begin{equation*}
		f_{yy}(y_j)
		\approx
		\frac{
			-f_{j+2}
			+16f_{j+1}
			-30f_j
			+16f_{j-1}
			-f_{j-2}
		}{
			12\Delta y^2
		}.
	\end{equation*}
	
	Near a Dirichlet boundary, one-sided polynomial finite-difference weights
	are used so that the spatial operator remains fourth-order accurate.
	High-order one-sided treatments are particularly relevant in
	electrochemical calculations because quantities such as current are often
	determined from concentration gradients at an electrode
	\cite{Bieniasz2003}. In the present implementation, six-point one-sided
	stencils are used at the first and last interior grid points, while the
	five-point centered stencil above is used in the
	interior.
	
	For fields satisfying a Neumann condition at the electrode, the boundary
	value is eliminated from a high-order one-sided approximation of the first
	derivative. Specifically, the relation
	\[
	f_y(0,\tau)=g(\tau)
	\]
	is imposed using the first six spatial nodes. The resulting expression for
	the boundary value is substituted into the fourth-order approximation of
	the Laplacian at the neighboring interior points. This implementation
	avoids a lower-order boundary closure and keeps the boundary
	discretization compatible with the fourth-order interior scheme. Related
	high-order boundary-gradient constructions for electrochemical
	reaction--diffusion problems are discussed in \cite{Bieniasz2003}.
	The same spatial operators are used for the HPM correction equations and
	for the independent nonlinear PDE benchmark, although the latter is solved
	on a substantially finer spatial grid.
	
	\subsection{Numerical Evaluation of Formulation I}
	\label{sec:num-form-I}
	
	The zeroth-order catalyst field $q_0$ is available analytically from
	\eqref{eq:I-q0closed}. The associated substrate field $a_0$ is reconstructed
	from the exact difference relation. Higher-order catalyst corrections satisfy
	the linear problems
	\begin{equation*}
		(q_n)_\tau
		=
		(q_n)_{yy}
		-\alpha q_n
		+
		F_n,
		\qquad
		F_n=\alpha q_{n-1}-R_{n-1},
		\qquad n\geq1,
	\end{equation*}
	with homogeneous initial and Dirichlet boundary conditions.
	
	For the production calculations, the full correction fields $q_n$ are
	computed from this correction equation using the fourth-order spatial
	method-of-lines discretization described above and adaptive
	\texttt{ode15s} integration. The Green-function representation derived in
	Section~\ref{sec:I-general} remains analytically exact and is used in two
	important ways: first, to validate the first correction independently; and
	second, to evaluate the surface derivative directly.
	
	\subsubsection{Structural reconstruction}
	
	For a sufficiently regular catalyst coefficient $f(y,\tau)$, define
	\begin{equation*}
		\mathcal H[f](y,\tau)
		=
		\int_0^\tau
		\frac{f_y(0,\sigma)}
		{\sqrt{\pi(\tau-\sigma)}}
		\exp\left[
		-\frac{y^2}{4(\tau-\sigma)}
		\right]
		\,d\sigma.
	\end{equation*}
	Then
	\[
	a_0=q_0+\gamma+\mathcal H[q_0],
	\]
	and, for $n\geq1$,
	\[
	a_n=q_n+\mathcal H[q_n].
	\]
	
	To remove the square-root endpoint factor, introduce
	\[
	\sigma=\tau\sin^2\theta,
	\qquad
	0\leq\theta\leq\frac{\pi}{2}.
	\]
	The reconstruction becomes
	\begin{equation}
		\mathcal H[f](y,\tau)
		=
		\frac{2\sqrt{\tau}}{\sqrt{\pi}}
		\int_0^{\pi/2}
		f_y\!\left(0,\tau\sin^2\theta\right)
		\sin\theta
		\exp\left[
		-\frac{y^2}{4\tau\cos^2\theta}
		\right]
		\,d\theta.
		\label{eq:num-H-transform}
	\end{equation}
	
	For the zeroth-order reference field, the surface derivative is known
	analytically:
	\begin{equation}
		(q_0)_y(0,\tau)
		=
		-m_0
		\left[
		\sqrt{\alpha}\,
		\erf\left(\sqrt{\alpha\tau}\right)
		+
		\frac{e^{-\alpha\tau}}{\sqrt{\pi\tau}}
		\right].
		\label{eq:num-q0y}
	\end{equation}
	Thus $a_0$ is reconstructed without numerical differentiation.
	
	For $n\geq1$, the boundary derivative is evaluated from the differentiated
	Green representation described below and then inserted into
	\eqref{eq:num-H-transform}. In the baseline calculations,
	$500$ quadrature points are used for the structural reconstruction.
	
	\subsubsection{Direct surface derivative}
	
	Although the full correction fields are computed by fourth-order MOL, the
	electrode derivative is evaluated directly from the Green representation.
	For
	\[
	F_n=\alpha q_{n-1}-R_{n-1},
	\]
	the regularized formula is
	\begin{equation}
		(q_n)_y(0,\tau)
		=
		\frac{4\sqrt{\tau}}{\sqrt{\pi}}
		\int_0^1
		e^{-\alpha\tau\zeta^2}
		\int_0^\infty
		z e^{-z^2}
		F_n\left(
		2\sqrt{\tau}\,\zeta z,
		\tau(1-\zeta^2)
		\right)
		\,dz\,d\zeta,
		\qquad n\geq1.
		\label{eq:num-current-I-direct}
	\end{equation}
	
	The rapidly decaying $z$-integral is truncated at
	\(
	z_{\max}=6.
	\)
	For the baseline calculations, Gauss--Legendre quadrature with
	\(
	N_\zeta=28,
	N_z=36
	\)
	is used. Consequently, the Formulation-I current is obtained from
	\begin{equation*}
		\psi_{\mathrm I}^{[N]}(\tau)
		=
		-
		\sum_{n=0}^{N}
		(q_n)_y(0,\tau),
	\end{equation*}
	with the zeroth-order derivative given analytically by
	\eqref{eq:num-q0y}.
	
	\subsubsection{Independent Green-function validation}
	\label{sec:num-I-validation}
	
	The first catalyst correction is also evaluated independently from its
	Green-function representation,
	\begin{equation*}
		q_1(y,\tau)
		=
		\int_0^\tau
		\int_0^\infty
		G_\alpha(y,\eta,\tau-\sigma)
		F_1(\eta,\sigma)
		\,d\eta\,d\sigma,
	\end{equation*}
	where
	\[
	F_1=\alpha q_0-R_0.
	\]
	After the regularizing transformation derived in
	Section~\ref{sec:I-first-order}, the time integral is evaluated by direct
	quadrature at selected values of $\tau$.
	
	The Green-function result is compared with the independently computed
	fourth-order MOL field. The corresponding substrate correction is obtained
	from the structural reconstruction and compared with a direct
	fourth-order coupled MOL solution of the equivalent frozen-coefficient
	system.
	
	Thus the validation compares
	\[
	q_1^{\mathrm{Green}}
	\quad\hbox{with}\quad
	q_1^{\mathrm{MOL4}},
	\]
	and
	\[
	a_1^{\mathrm{Green/reconstruction}}
	\quad\hbox{with}\quad
	a_1^{\mathrm{MOL4,coupled}}.
	\]
	This comparison is independent of the production field calculation and
	therefore provides a useful check of both the Green representation and the
	high-order MOL implementation.
	
	\subsection{Numerical Evaluation of Formulation II}
	\label{sec:num-form-II}
	
	For Formulation~II, the reference quantities
	\[
	\bar q,\qquad
	\bar a,\qquad
	\bar D,\qquad
	\bar R,\qquad
	\beta_a,\qquad
	\beta_q,\qquad
	\mathcal E
	\]
	are evaluated first on the space--time grid.
	
	Define
	\[
	H_0=-\mathcal E,
	\]
	and, for $n\geq1$,
	\[
	H_n=-\mathcal N_{n-1}.
	\]
	Then each correction pair satisfies
	\begin{equation}
		\begin{aligned}
			(u_n)_\tau
			&=
			(u_n)_{yy}
			-\beta_q u_n
			-\beta_a v_n
			+
			H_n,
			\\
			(v_n)_\tau
			&=
			(v_n)_{yy}
			-\beta_q u_n
			-\beta_a v_n
			+
			H_n,
		\end{aligned}
		\qquad n\geq0.
		\label{eq:num-II-general}
	\end{equation}
	
	The corrections satisfy
	\[
	u_n(y,0)=0,
	\qquad
	v_n(y,0)=0,
	\]
	with
	\[
	u_n(0,\tau)=0,
	\qquad
	u_n(L,\tau)=0,
	\]
	and
	\[
	(v_n)_y(0,\tau)=0,
	\qquad
	v_n(L,\tau)=0.
	\]
	
	The two fields are discretized simultaneously using the fourth-order
	spatial operators described in Section~\ref{sec:num-spatial}. Only interior
	unknowns are advanced. The Dirichlet values of $u_n$ are imposed directly,
	whereas the electrode value of $v_n$ is eliminated through the high-order
	Neumann relation.
	
	If
	\[
	X_n(\tau)
	=
	\begin{bmatrix}
		U_n(\tau)\\
		V_n(\tau)
	\end{bmatrix},
	\]
	denotes the vector of semidiscrete correction unknowns, the resulting
	system may be written schematically as
	\begin{equation*}
		\frac{dX_n}{d\tau}
		=
		M(\tau)X_n+\mathbf H_n(\tau),
	\end{equation*}
	where $M(\tau)$ contains the fourth-order diffusion matrices together with
	the time-dependent multiplication operators associated with
	$-\beta_q$ and $-\beta_a$.
	
	This semidiscrete system is integrated by
	\texttt{ode15s}, with
	\[
	\mathrm{RelTol}=10^{-9},
	\qquad
	\mathrm{AbsTol}=10^{-11}.
	\]
	The prescribed output-time spacing is also used as the maximum internal
	time step.
	
	After $(u_n,v_n)$ has been computed, the nonlinear coefficient
	$\mathcal N_n$ is evaluated from the recurrence in
	Section~\ref{sec:II-nonlinear-recursion}. The calculation therefore proceeds
	as
	\[
	(\bar q,\bar a)
	\longrightarrow
	(\beta_a,\beta_q,\mathcal E)
	\longrightarrow
	(u_0,v_0)
	\longrightarrow
	\mathcal N_0
	\longrightarrow
	(u_1,v_1)
	\longrightarrow
	\mathcal N_1
	\longrightarrow
	\cdots .
	\]
	
	\subsubsection{Direct heat-kernel surface derivative}
	
	The Formulation-II current is not obtained by differentiating the discrete
	correction fields. From \eqref{eq:num-II-general},
	\begin{equation*}
		(u_n)_\tau-(u_n)_{yy}
		=
		\mathcal S_n(y,\tau),
	\end{equation*}
	where
	\begin{equation*}
		\mathcal S_0
		=
		-\beta_q u_0-\beta_a v_0-\mathcal E,
	\end{equation*}
	and, for $n\geq1$,
	\begin{equation*}
		\mathcal S_n
		=
		-\beta_q u_n-\beta_a v_n-\mathcal N_{n-1}.
	\end{equation*}
	
	Differentiating the half-line Dirichlet heat-kernel representation gives
	\begin{equation*}
		(u_n)_y(0,\tau)
		=
		\frac{4\sqrt{\tau}}{\sqrt{\pi}}
		\int_0^1
		\int_0^\infty
		z e^{-z^2}
		\mathcal S_n
		\left(
		2\sqrt{\tau}\,\zeta z,
		\tau(1-\zeta^2)
		\right)
		\,dz\,d\zeta .
	\end{equation*}
	
	For the baseline calculations, the $z$-integral is truncated at
	\(
	z_{\max}=7,
	\)
	and the quadrature uses
	\(
	N_\zeta=48,
	N_z=64.
	\)
	Since
	\(
	\bar q=q_0=q_{\mathrm P},
	\)
	the reference contribution is known analytically:
	\begin{equation*}
		-\bar q_y(0,\tau)
		=
		m_0
		\left[
		\sqrt{\alpha}\,
		\erf\left(\sqrt{\alpha\tau}\right)
		+
		\frac{e^{-\alpha\tau}}
		{\sqrt{\pi\tau}}
		\right].
	\end{equation*}
	Therefore,
	\begin{equation}
		\psi_{\mathrm{II}}^{[N]}(\tau)
		=
		-\bar q_y(0,\tau)
		-
		\sum_{n=0}^{N-1}
		(u_n)_y(0,\tau),
		\qquad
		N\geq1.
		\label{eq:num-II-current-direct}
	\end{equation}
	
	A fourth-order one-sided derivative of the discrete Formulation-II field is
	also computed, but only as a diagnostic. All Formulation-II current errors
	reported in Section~\ref{sec:numerical-results} use the direct heat-kernel
	formula \eqref{eq:num-II-current-direct}.
	
	\subsubsection{Structural consistency check}
	
	For each correction level, the directly solved substrate correction is
	compared with the structural reconstruction
	\begin{equation*}
		v_n^{\mathrm{rec}}(y,\tau)
		=
		u_n(y,\tau)
		+
		\mathcal H[u_n](y,\tau),
		\qquad n\geq0.
	\end{equation*}
	The reconstruction uses \eqref{eq:num-H-transform}. The coefficient-level
	consistency error is measured over the post-transient interval by
	\begin{equation*}
		E_{\mathrm{rec}}^{(n)}
		=
		\max_{\tau\geq\tau_c}
		\max_{0\leq y\leq L}
		\left|
		v_n^{\mathrm{dir}}(y,\tau)
		-
		v_n^{\mathrm{rec}}(y,\tau)
		\right|.
	\end{equation*}
	
	For Formulation~I, an analogous independent check is performed by comparing
	the structurally reconstructed substrate coefficients with those obtained
	from the equivalent coupled frozen-coefficient MOL system. These checks
	test the numerical implementation of the exact difference structure and
	are distinct from comparison with the nonlinear PDE benchmark.
	
	\subsection{Independent Nonlinear PDE Reference Solution}
	\label{sec:nonlinear-reference}
	
	The published and HPM approximations are assessed against an independent
	solution of the full nonlinear reduced system. Independent numerical
	benchmarks are important for nonlinear electroanalytical
	reaction--diffusion models because agreement of an approximation with its
	own auxiliary or linearized problem does not establish accuracy for the
	original nonlinear PDE
	\cite{VynnyckyMcKeeBieniasz2021}.
	We therefore solve the full nonlinear reduced system
	\begin{equation*}
		\begin{aligned}
			q_\tau
			&=
			q_{yy}
			-
			\frac{\lambda_1aq}{1+\rho(1-q)},
			\\
			a_\tau
			&=
			a_{yy}
			-
			\frac{\lambda_1aq}{1+\rho(1-q)}.
		\end{aligned}
	\end{equation*}
	
	The truncated-domain conditions are
	\[
	q(0,\tau)=m_0,
	\qquad
	q(L,\tau)=0,
	\]
	and
	\[
	a_y(0,\tau)=0,
	\qquad
	a(L,\tau)=\gamma.
	\]
	
	The benchmark uses the same fourth-order spatial discretization as the HPM
	calculations but on an independently refined spatial grid. The resulting
	nonlinear semidiscrete system is integrated using \texttt{ode15s}. The
	nonlinear reaction term is evaluated pointwise from the current numerical
	values of $a$ and $q$.
	
	For the baseline calculations,
	\[
	N_y^{\mathrm{ref}}=1281,
	\]
	and
	\[
	\mathrm{RelTol}=10^{-9},
	\qquad
	\mathrm{AbsTol}=10^{-11}.
	\]
	The fine-grid solution is interpolated onto the HPM spatial grid only after
	the nonlinear PDE has been solved. The reference electrode current is
	differentiated on the fine grid before any interpolation.
	
	The initial condition
	\(
	q(y,0)=0
	\)
	and the boundary condition
	\(
	q(0,\tau)=m_0,
	\tau>0,
	\)
	remain incompatible at the corner $(0,0)$, as in the continuous problem.
	No artificial smoothing of this corner is introduced.
	
	\subsection{Evaluation of the Electrode Current}
	\label{sec:num-current}
	
	The dimensionless current is
	\begin{equation*}
		\psi(\tau)=-q_y(0,\tau).
	\end{equation*}
	Because the current is determined by a boundary derivative, its numerical
	evaluation is generally more sensitive to spatial resolution and boundary
	discretization than the concentration field itself. High-order one-sided
	approximations to electrode concentration gradients have therefore been
	studied specifically for electrochemical reaction--diffusion simulations
	\cite{Bieniasz2003}.
	
	For the published approximation,
	\begin{equation*}
		\psi_{\mathrm P}(\tau)
		=
		m_0
		\left[
		\sqrt K\,
		\erf\left(\sqrt{K\tau}\right)
		+
		\frac{e^{-K\tau}}{\sqrt{\pi\tau}}
		\right].
	\end{equation*}
	
	For Formulation~I, the zeroth-order term is analytic and every higher-order
	surface derivative is evaluated from \eqref{eq:num-current-I-direct}. For
	Formulation~II, all reported currents use the direct heat-kernel formula
	\eqref{eq:num-II-current-direct}.
	
	For the nonlinear benchmark, the catalyst field is available on the fine
	reference grid. Its surface derivative is evaluated using the fourth-order
	one-sided formula
	\begin{equation*}
		q_y(0,\tau)
		\approx
		\frac{
			-25Q_0
			+48Q_1
			-36Q_2
			+16Q_3
			-3Q_4
		}{
			12\Delta y_{\mathrm{ref}}
		},
	\end{equation*}
	where
	\[
	Q_j(\tau)=q(y_j,\tau).
	\]
	Thus
	\begin{equation*}
		\psi^{\mathrm{PDE}}(\tau)
		\approx
		\frac{
			25Q_0
			-48Q_1
			+36Q_2
			-16Q_3
			+3Q_4
		}{
			12\Delta y_{\mathrm{ref}}
		}.
	\end{equation*}
	
	\subsection{Error Measures}
	\label{sec:num-errors}
	
	Let
	\[
	q^{\mathrm{PDE}}(y,\tau),
	\qquad
	a^{\mathrm{PDE}}(y,\tau),
	\qquad
	\psi^{\mathrm{PDE}}(\tau)
	\]
	denote the independent nonlinear reference solution.
	
	For an approximation
	\[
	(q^{\mathrm{app}},a^{\mathrm{app}},\psi^{\mathrm{app}}),
	\]
	define the time-dependent concentration errors
	\begin{equation*}
		E_q(\tau)
		=
		\max_{0\leq y\leq L}
		\left|
		q^{\mathrm{app}}(y,\tau)
		-
		q^{\mathrm{PDE}}(y,\tau)
		\right|,
	\end{equation*}
	and
	\begin{equation*}
		E_a(\tau)
		=
		\max_{0\leq y\leq L}
		\left|
		a^{\mathrm{app}}(y,\tau)
		-
		a^{\mathrm{PDE}}(y,\tau)
		\right|.
	\end{equation*}
	The current error is
	\begin{equation*}
		E_\psi(\tau)
		=
		\left|
		\psi^{\mathrm{app}}(\tau)
		-
		\psi^{\mathrm{PDE}}(\tau)
		\right|.
	\end{equation*}
	
	Because
	\[
	q_y(0,\tau)=O(\tau^{-1/2})
	\qquad
	(\tau\to0^+),
	\]
	the very short-time current error is particularly sensitive to numerical
	resolution of the initial corner layer. The principal scalar errors are
	therefore measured over
	\[
	\tau_c\leq\tau\leq\tau_{\max},
	\]
	where
	\begin{equation*}
		\tau_c=0.05\,\tau_{\max}.
	\end{equation*}
	
	The reported errors are
	\begin{equation*}
		E_q
		=
		\max_{\tau_c\leq\tau\leq\tau_{\max}}
		\max_{0\leq y\leq L}
		\left|
		q^{\mathrm{app}}(y,\tau)
		-
		q^{\mathrm{PDE}}(y,\tau)
		\right|,
	\end{equation*}
	\begin{equation*}
		E_a
		=
		\max_{\tau_c\leq\tau\leq\tau_{\max}}
		\max_{0\leq y\leq L}
		\left|
		a^{\mathrm{app}}(y,\tau)
		-
		a^{\mathrm{PDE}}(y,\tau)
		\right|,
	\end{equation*}
	and
	\begin{equation*}
		E_\psi
		=
		\max_{\tau_c\leq\tau\leq\tau_{\max}}
		\left|
		\psi^{\mathrm{app}}(\tau)
		-
		\psi^{\mathrm{PDE}}(\tau)
		\right|.
	\end{equation*}
	
	The full positive-time histories are retained separately so that the
	corner-layer behavior can be distinguished from the post-transient
	approximation error.
	
	\subsection{Baseline Numerical Parameters and Convergence Checks}
	\label{sec:num-parameters}
	
	Unless otherwise stated, the baseline physical parameters are
	\[
	\lambda_1=1,
	\qquad
	\rho=1,
	\qquad
	\gamma=1,
	\qquad
	\xi=0,
	\]
	and hence
	\[
	m_0=\frac12.
	\]
	
	The baseline numerical parameters are summarized in
	Table~\ref{tab:numerical-parameters}.
	
	\begin{table}[htbp]
		\centering
		\small
		\renewcommand{\arraystretch}{1.25}
		\begin{tabular}{ll}
			\toprule
			\textbf{Quantity} & \textbf{Baseline value}\\
			\midrule
			Final time $\tau_{\max}$ & $1$\\
			Domain length $L$ & $8$\\
			HPM spatial grid points $N_y$ & $161$\\
			Output time points $N_t$ & $401$\\
			HPM truncation label $N$ & $3$\\
			Reference spatial grid points $N_y^{\mathrm{ref}}$ & $1281$\\
			Structural-reconstruction quadrature points & $500$\\
			Formulation-I surface quadrature $(N_\zeta,N_z)$ & $(28,36)$\\
			Formulation-I surface cutoff $z_{\max}$ & $6$\\
			Formulation-II surface quadrature $(N_\zeta,N_z)$ & $(48,64)$\\
			Formulation-II surface cutoff $z_{\max}$ & $7$\\
			Linear MOL relative tolerance & $10^{-9}$\\
			Linear MOL absolute tolerance & $10^{-11}$\\
			Nonlinear reference relative tolerance & $10^{-9}$\\
			Nonlinear reference absolute tolerance & $10^{-11}$\\
			Post-transient start $\tau_c$ & $0.05$\\
			\bottomrule
		\end{tabular}
		\caption{Baseline numerical parameters used in the comparisons.}
		\label{tab:numerical-parameters}
	\end{table}
	
	Separate numerical studies are used to assess domain truncation, spatial
	resolution, temporal output resolution, and structural-reconstruction
	quadrature. Specifically,
	\[
	L=6,\;8,\;10,\;12,
	\]
	\[
	N_y=81,\;121,\;161,\;241,\;321,
	\]
	\[
	N_t=201,\;401,\;801,
	\]
	and
	\[
	N_{\mathrm{quad}}=200,\;500,\;1000
	\]
	are examined independently.
	A further high-resolution spatial study uses
	\[
	N_y=161,\;241,\;321,\;481,\;641
	\]
	while retaining the same independently computed fine nonlinear benchmark.
	This additional experiment is used to distinguish HPM truncation error from
	the numerical error floor in the concentration and current calculations.
	The HPM truncation order is varied independently:
	\(
	N=0,\;1,\;2,\;3.
	\)
	Thus variation with $N$ measures convergence of the analytical approximation,
	whereas variation with $L$, $N_y$, $N_t$, and quadrature resolution measures
	numerical evaluation sensitivity.
	
	\subsection{Computational Efficiency Measures}
	\label{sec:num-efficiency-method}
	
	Runtime alone is not sufficient for comparing the three numerical
	approaches because it depends on hardware, operating-system load, and solver
	implementation. The efficiency study therefore reports both repeated timing
	measurements and hardware-independent solver-work quantities.
	
	For the two HPM methods, the comparison is made at $N=1$, where the
	low-order advantage of the profile-linearized formulation is most relevant.
	Three dedicated runs are performed for each method, and the median, minimum,
	and maximum elapsed runtime are recorded.
	
	The semidiscrete ODE state dimensions are
	\[
	N_{\mathrm{state}}^{\mathrm I}
	=
	N_y-2,
	\]
	\[
	N_{\mathrm{state}}^{\mathrm{II}}
	=
	2(N_y-2),
	\]
	and
	\[
	N_{\mathrm{state}}^{\mathrm{PDE}}
	=
	2(N_y^{\mathrm{ref}}-2).
	\]
	
	In addition, the actual number of right-hand-side evaluations performed by
	\texttt{ode15s}, denoted $N_{\mathrm{RHS}}$, is counted directly. We define
	the solver-work index
	\begin{equation*}
		W
		=
		N_{\mathrm{state}}N_{\mathrm{RHS}}.
	\end{equation*}
	This is not a floating-point operation count, but it provides a reproducible
	measure of the size and amount of ODE work that is substantially less
	hardware dependent than runtime.
	Kernel-based postprocessing is reported separately. For Formulation~I, each
	positive output time requires
	\(
	28\times36=1008
	\)
	surface-flux quadrature nodes per retained correction. For
	Formulation~II, the corresponding number is
	\(
	48\times64=3072.
	\)
	The structural reconstruction uses $500$ quadrature nodes per positive
	output time for each reconstructed coefficient. These quadrature operations
	are not included in the work index $W$ and are
	therefore reported separately in the efficiency discussion.
	
	The accuracy, convergence, structural checks, parameter studies, and
	efficiency results are presented in
	Section~\ref{sec:numerical-results}.
	
	\section{Numerical Results}
	\label{sec:numerical-results}
	
	This section evaluates the published low-scan-rate approximation,
	Formulation~I, and Formulation~II against an independent numerical
	solution of the full nonlinear reaction--diffusion system.  The numerical
	methodology used to compute the HPM corrections and the nonlinear benchmark
	was described in the preceding section.  Particular attention is given to
	four questions:
	\begin{enumerate}\itemsep -0.1cm
		\item whether the Green-function and high-order MOL implementations are
		mutually consistent;
		\item how rapidly the two HPM formulations converge with perturbation
		order;
		\item how robust the approximations remain as the reaction parameters
		$\lambda_1$ and $\rho$ vary; and
		\item whether the improved accuracy is obtained at substantially lower
		computational cost than direct solution of the nonlinear PDE.
	\end{enumerate}
	
	\subsection{Baseline parameters and error measures}
	\label{subsec:baseline-parameters}
	
	Unless otherwise stated, the baseline parameters are
	\[
	\gamma=1,\qquad
	m_0=\frac12,\qquad
	\lambda_1=1,\qquad
	\rho=1,
	\]
	with
	\[
	L=8,\qquad
	0\leq \tau\leq 1.
	\]
	The HPM calculations use
	\[
	N_y=161,\qquad
	N_\tau=401,
	\]
	while the independent nonlinear PDE benchmark is computed on the finer grid
	\[
	N_y^{\mathrm{ref}}=1281.
	\]
	Both the linear correction equations and the nonlinear reference problem are
	integrated using the fourth-order spatial MOL discretization described
	previously, together with adaptive \texttt{ode15s} time integration.  The
	relative and absolute tolerances are
	\[
	\mathrm{RelTol}=10^{-9},
	\qquad
	\mathrm{AbsTol}=10^{-11}.
	\]
	The structural reconstruction uses $N_{\mathrm{quad}}=500$ quadrature
	points.
	
	Because the initial and electrode data generate a sharp corner layer near
	$\tau=0$, quantitative errors are measured over the post-transient interval
	\[
	\tau\geq\tau_c,\qquad \tau_c=0.05.
	\]
	For an approximation $(q^{\mathrm{app}},a^{\mathrm{app}},\psi^{\mathrm{app}})$,
	we define
	\begin{align*}
		E_q
		&=
		\max_{\tau\geq\tau_c}
		\max_{0\leq y\leq L}
		\left|
		q^{\mathrm{app}}(y,\tau)-q^{\mathrm{PDE}}(y,\tau)
		\right|,
		\\
		E_a
		&=
		\max_{\tau\geq\tau_c}
		\max_{0\leq y\leq L}
		\left|
		a^{\mathrm{app}}(y,\tau)-a^{\mathrm{PDE}}(y,\tau)
		\right|,
		\\
		E_\psi
		&=
		\max_{\tau\geq\tau_c}
		\left|
		\psi^{\mathrm{app}}(\tau)-\psi^{\mathrm{PDE}}(\tau)
		\right|.
	\end{align*}
	
	\subsection{Independent validation of Formulation I}
	\label{subsec:green-validation}
	
	Before comparing the HPM approximations with the nonlinear PDE, we first
	verify the numerical realization of Formulation~I.  The first correction
	$q_1$ is computed independently in two ways: directly from its
	Green-function representation and from the fourth-order MOL discretization.
	The corresponding substrate correction $a_1$ is obtained from the structural
	reconstruction in the Green-function calculation and from the equivalent
	coupled frozen-coefficient MOL problem in the independent calculation.
	
\begin{figure}[H]
		\centering
		\includegraphics[width=0.9\linewidth, trim=0 180 0 180, clip]
		{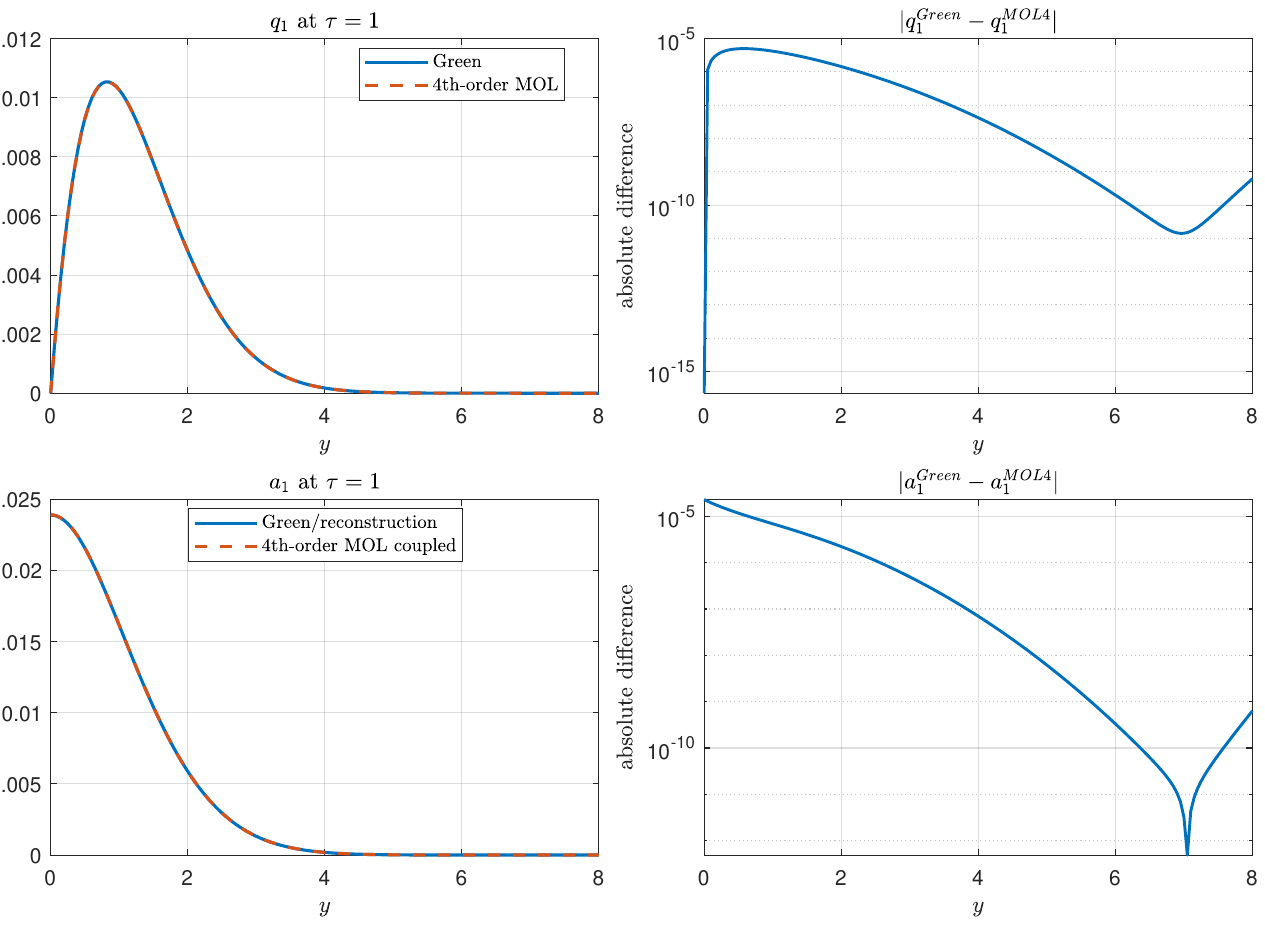}
		\caption{Independent validation of the first Formulation~I correction.
			The left panels compare the Green-function/reconstruction solution with
			the fourth-order MOL solution at $\tau=1$, while the right panels show
			the corresponding absolute differences.}
		\label{fig:green-mol4-validation}
	\end{figure}

	Figure~\ref{fig:green-mol4-validation} shows the comparison at $\tau=1$.
	The two profiles are visually indistinguishable.  Over the validation times
	\[
	\tau=0.10,\;0.25,\;0.50,\;1.00,
	\]
	the largest discrepancies are approximately
	\[
	\max |q_1^{\mathrm{Green}}-q_1^{\mathrm{MOL4}}|
	\leq 5.03\times10^{-6},
	\]
	and
	\[
	\max |a_1^{\mathrm{Green}}-a_1^{\mathrm{MOL4}}|
	\leq 2.35\times10^{-5}.
	\]
	This agreement provides an independent check of both the Green-function
	representation and the fourth-order MOL implementation.
	
	\subsection{Representative concentration and current profiles}
	\label{subsec:representative-results}
	
\begin{figure}[htbp]
		\centering
		\includegraphics[width=0.9\linewidth, trim=0 260 0 260, clip ]
		{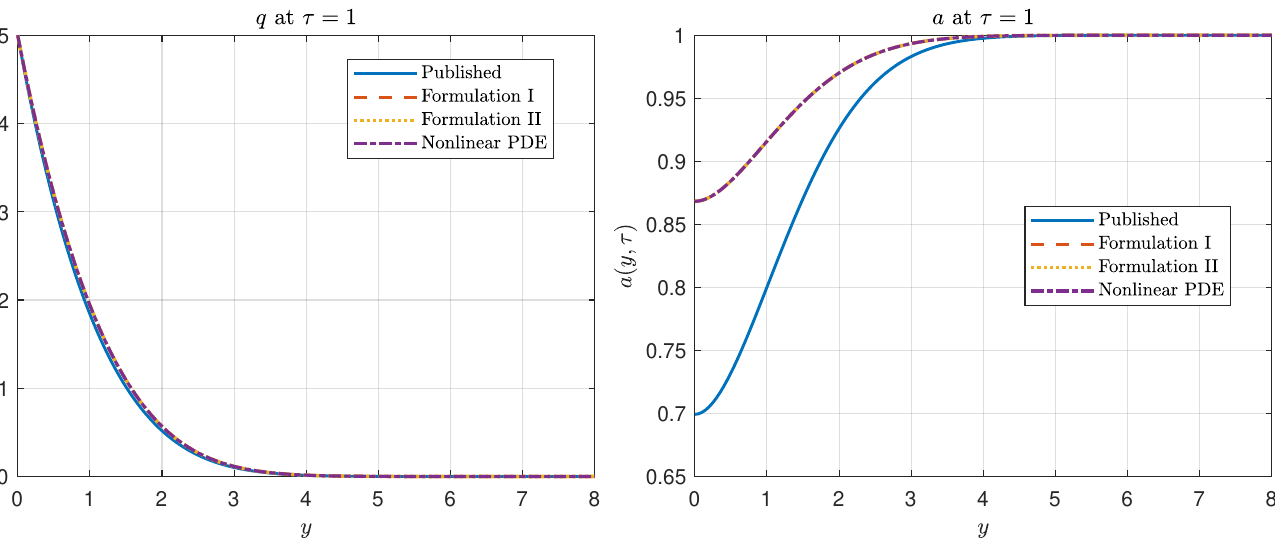}
		\caption{Representative catalyst and substrate profiles at $\tau=1$
			for the published approximation, Formulation~I, Formulation~II, and the
			nonlinear PDE benchmark.}
		\label{fig:representative-concentrations}
	\end{figure}

	Figure~\ref{fig:representative-concentrations} compares the published
	approximation, the two proposed HPM formulations, and the nonlinear PDE
	solution at $\tau=1$.  For the catalyst concentration $q$, all four curves
	have the same overall diffusion-layer structure, although the discrepancy of
	the published approximation is visible upon quantitative comparison.  The
	difference is substantially more pronounced for the substrate concentration
	$a$.  In particular, the published approximation underestimates the
	substrate concentration near the electrode, whereas both structure-preserving
	HPM approximations are nearly indistinguishable from the nonlinear PDE
	solution.
	
	For the baseline parameters, the maximum post-transient errors are summarized
	in Table~\ref{tab:baseline-errors}.  The published substrate approximation has
	an error of approximately $1.69\times10^{-1}$, whereas the errors from
	Formulations~I and II are reduced to approximately $2.09\times10^{-5}$ and
	$3.61\times10^{-6}$, respectively.
	
	\begin{table}[htbp]
		\centering
		\caption{Baseline maximum errors over $\tau\geq0.05$.
			The HPM results use $N=3$.}
		\label{tab:baseline-errors}
		\small
		\begin{tabular}{lccc}
			\hline
			Method & $E_q$ & $E_a$ & $E_\psi$\\
			\hline
			Published
			& $1.0240\times10^{-2}$
			& $1.6896\times10^{-1}$
			& $2.8073\times10^{-2}$\\
			Formulation I
			& $5.1444\times10^{-7}$
			& $2.0865\times10^{-5}$
			& $5.3916\times10^{-5}$\\
			Formulation II
			& $1.5508\times10^{-7}$
			& $3.6082\times10^{-6}$
			& $5.3489\times10^{-5}$\\
			\hline
		\end{tabular}
	\end{table}
	
\begin{figure}[htbp]
		\centering
		\includegraphics[width=0.9\linewidth, trim=0 280 0 280, clip]
		{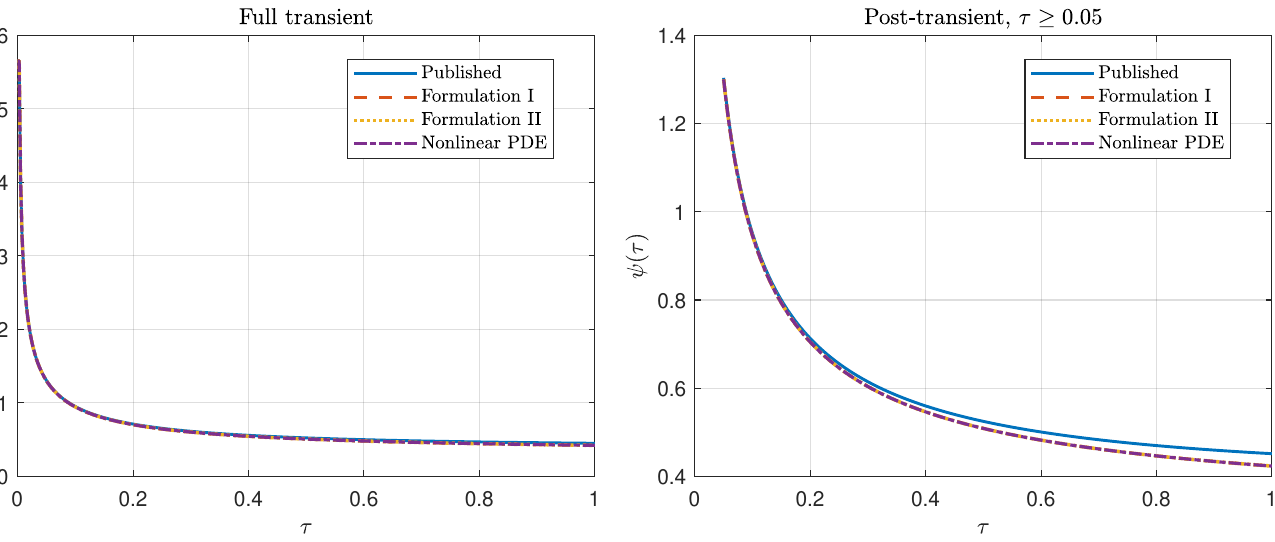}
		\caption{Electrode current for the published approximation,
			Formulation~I, Formulation~II, and the nonlinear PDE benchmark.
			The right panel enlarges the post-transient interval
			$\tau\geq0.05$.}
		\label{fig:representative-current}
	\end{figure}

	The corresponding current curves are shown in
	Figure~\ref{fig:representative-current}.  All methods reproduce the dominant
	short-time Cottrell-type behavior, but the difference between the published
	current and the nonlinear PDE becomes increasingly visible after the initial
	transient.  Both HPM currents closely follow the nonlinear reference.
	
\begin{figure}[htbp]
		\centering
		\includegraphics[width=0.9\linewidth, trim=0 260 0 260, clip]
		{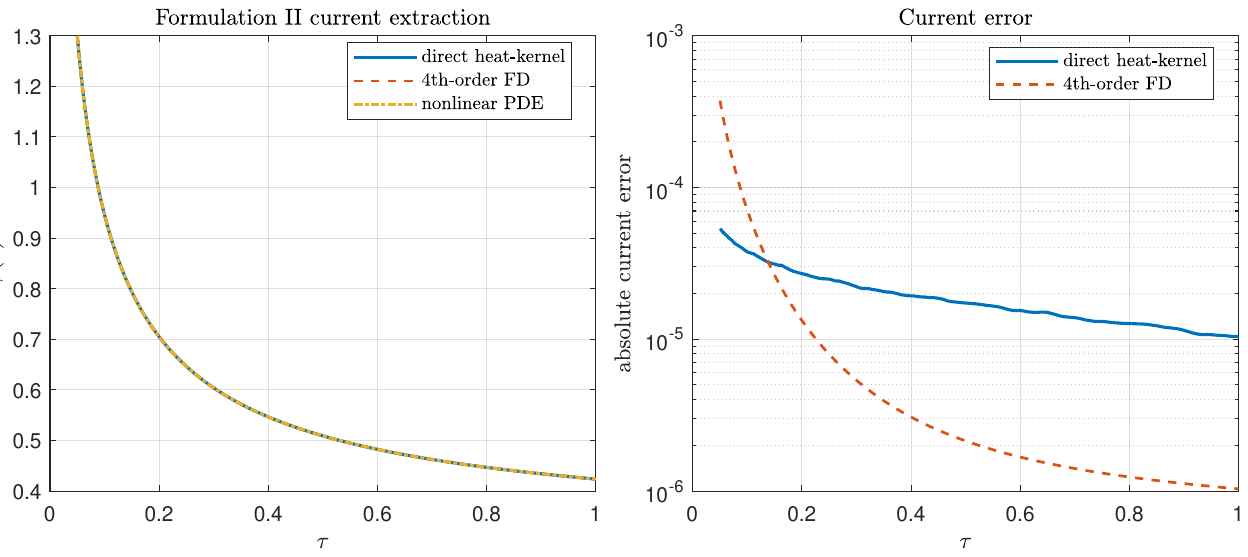}
		\caption{Comparison of direct heat-kernel and fourth-order
			finite-difference current extraction for Formulation~II.
			The nonlinear PDE current is included as the reference.}
		\label{fig:current-extraction}
	\end{figure}

	For Formulation~II, the reported current is evaluated directly from the
	differentiated heat-kernel representation rather than by differentiating the
	discrete concentration field.  As shown in
	Figure~\ref{fig:current-extraction}, this direct evaluation is particularly
	advantageous during the short-time boundary-layer regime.  Over
	$\tau\geq0.05$, the maximum direct-kernel current error is approximately
	$5.35\times10^{-5}$, compared with approximately $3.73\times10^{-4}$ for the
	fourth-order finite-difference extraction.  At later times, after the
	solution becomes smooth, the finite-difference derivative becomes extremely
	accurate.  The direct heat-kernel result is nevertheless retained as the
	reported Formulation~II current because it follows directly from the
	continuous correction equation and avoids numerical differentiation at the
	electrode.
	
	\subsection{Convergence with HPM order}
	\label{subsec:hpm-order}
	
\begin{figure}[htbp]
		\centering
		\includegraphics[width=0.9\linewidth, trim=0 260 0 260, clip]
		{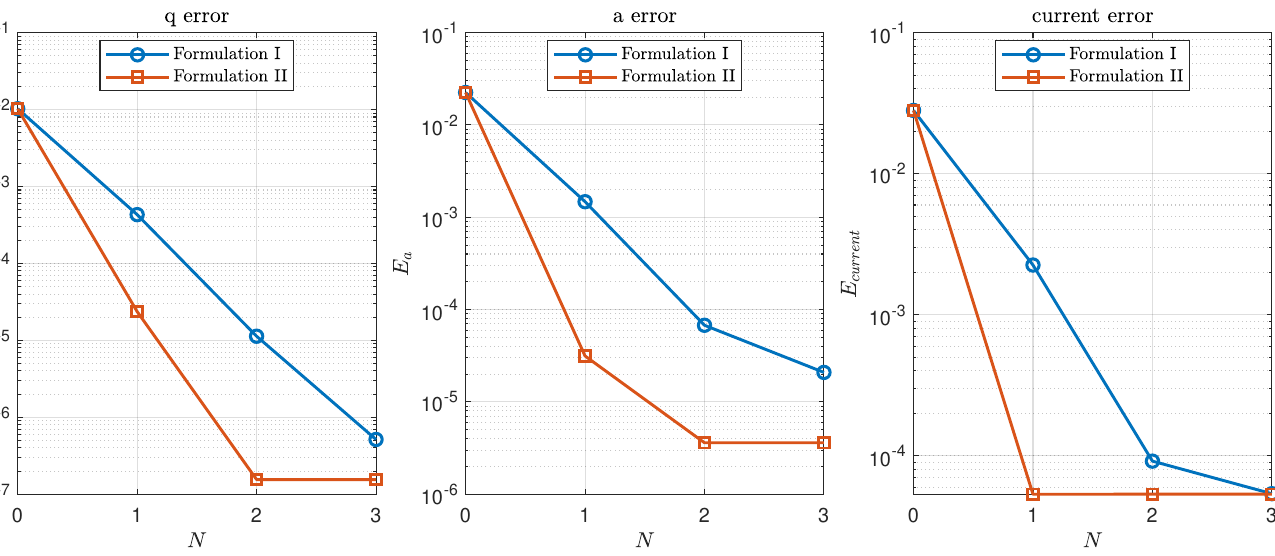}
		\caption{Maximum errors as functions of HPM truncation order.
			Formulation~II reaches the numerical error floor at substantially lower
			order than Formulation~I.}
		\label{fig:hpm-order}
	\end{figure}

	Figure~\ref{fig:hpm-order} shows the error as the HPM truncation order is
	increased from $N=0$ to $N=3$.  Formulation~I exhibits progressive reduction
	of all three errors.  Formulation~II converges considerably more rapidly:
	already at $N=1$,
	\[
	E_q^{II}=2.36\times10^{-5},
	\qquad
	E_a^{II}=3.14\times10^{-5},
	\qquad
	E_\psi^{II}=5.35\times10^{-5}.
	\]
	At the same order, Formulation~I gives
	\[
	E_q^{I}=4.31\times10^{-4},
	\qquad
	E_a^{I}=1.47\times10^{-3},
	\qquad
	E_\psi^{I}=2.25\times10^{-3}.
	\]
	Thus, at $N=1$, Formulation~II improves the catalyst, substrate, and current
	errors by factors of approximately
	\(
	18.2, 46.9, 42.1,
	\)
	respectively.
	
	A distinction at $N=0$ is important.  Both formulations use
	\(
	q_0=q_P,
	\psi_0=\psi_P,
	\)
	but the corresponding substrate field is the structure-preserving
	reconstruction
	\(
	a_0=q_0+\gamma+\mathcal H[q_0],
	\)
	rather than the separately approximated published substrate $a_P$.
	Consequently, the $N=0$ substrate error in
	Figure~\ref{fig:hpm-order},
	\[
	E_a^{[0]}\approx2.24\times10^{-2},
	\]
	is already much smaller than the published substrate error
	$1.69\times10^{-1}$.
	
	For Formulation~II, the concentration errors essentially stop decreasing
	after $N=2$.  This behavior is consistent with the remaining error approaching a combined numerical floor associated with spatial discretization, interpolation, quadrature, and the finite-resolution reference solution. The higher-order spatial study below confirms this interpretation.
	
	\subsection{Numerical sensitivity and spatial refinement}
	\label{subsec:numerical-convergence}
	
	The domain length, structural quadrature, and temporal resolution were varied
	independently before fixing the baseline numerical parameters.  Increasing
	$L$ from $6$ to $8$ substantially reduces the catalyst truncation error,
	whereas the results for
	\[
	L=8,\;10,\;12
	\]
	are essentially unchanged.  We therefore use $L=8$.
	Similarly, increasing the structural quadrature from
	$N_{\mathrm{quad}}=200$ to $500$ reduces the Formulation~II substrate
	reconstruction error, while the change from $500$ to $1000$ is negligible.
	The choice
	\(
	N_{\mathrm{quad}}=500
	\)
	is therefore sufficient.
	
	Temporal refinement from $N_\tau=201$ to $401$ improves the principal error
	measures.  Further refinement to $N_\tau=801$ produces only small,
	non-monotone changes once the errors have reached the combined spatial,
	interpolation, quadrature, and benchmark error floor.  For this reason,
	$N_\tau=401$ is used in the broad parameter studies.
	
	\begin{figure}[p]
		\centering
		\includegraphics[width=0.9\linewidth, trim=0 260 0 260, clip]
		{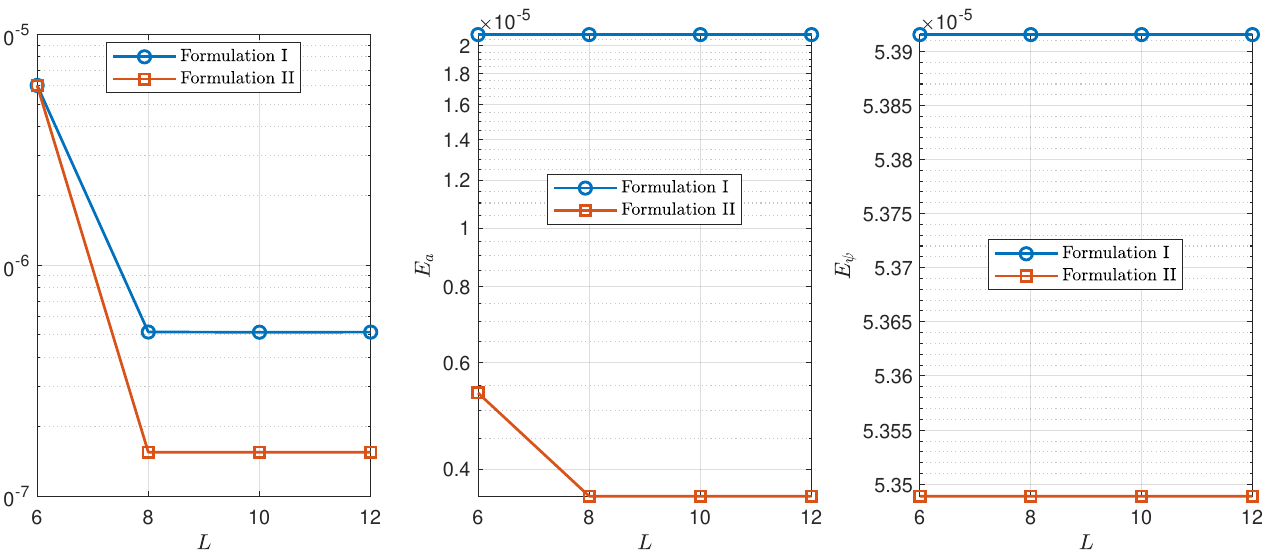}
		
		\vspace{0.4cm}
		
		\includegraphics[width=0.9\linewidth, trim=0 260 0 260, clip ]
		{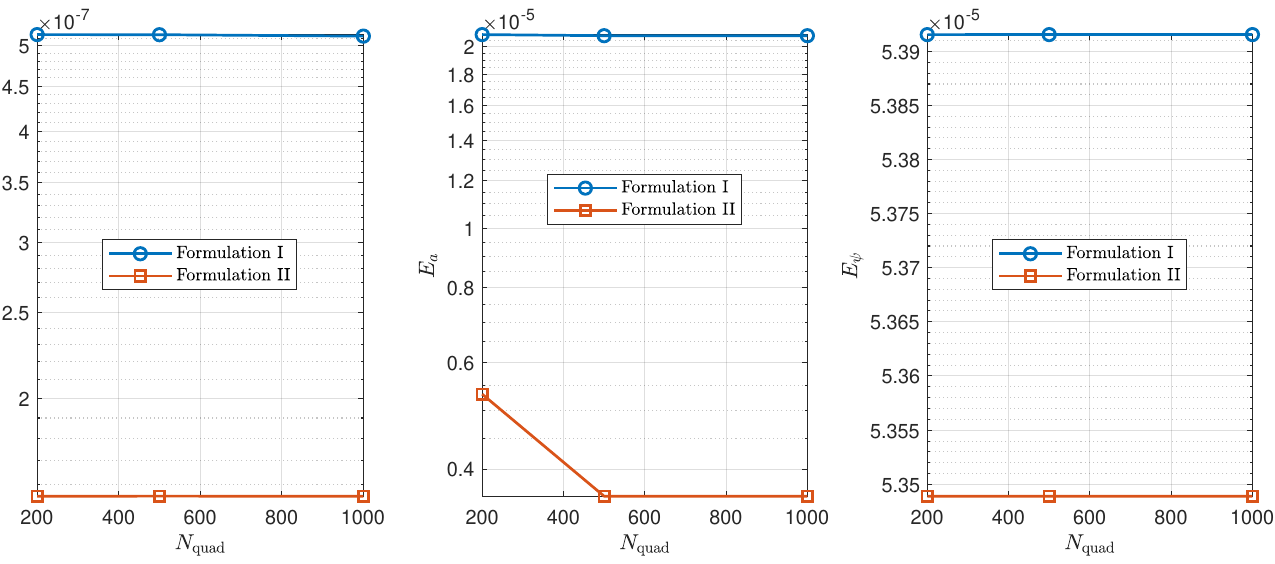}
		
		\vspace{0.4cm}
		
		\includegraphics[width=0.9\linewidth, trim=0 260 0 260, clip ]
		{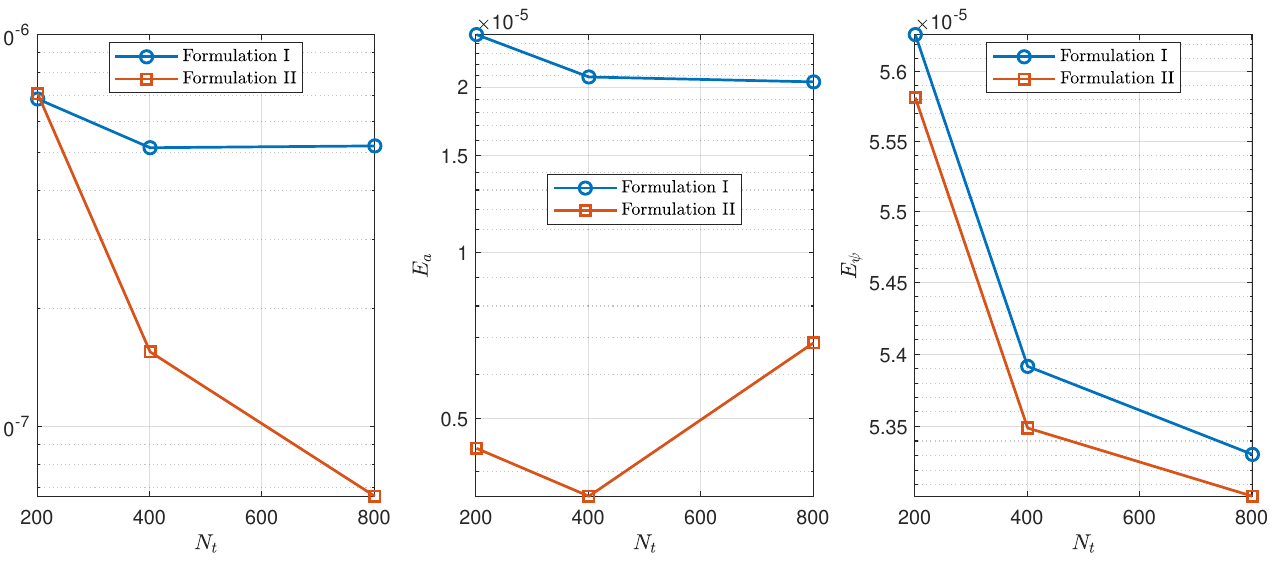}
		\caption{Sensitivity of the numerical errors to domain truncation
			$L$ (top), structural quadrature size $N_{\mathrm{quad}}$ (middle),
			and temporal output resolution $N_\tau$ (bottom).}
		\label{fig:basic-convergence}
	\end{figure}

\begin{figure}[htbp]
		\centering
		\includegraphics[width=0.9\linewidth, trim=0 260 0 260, clip]
		{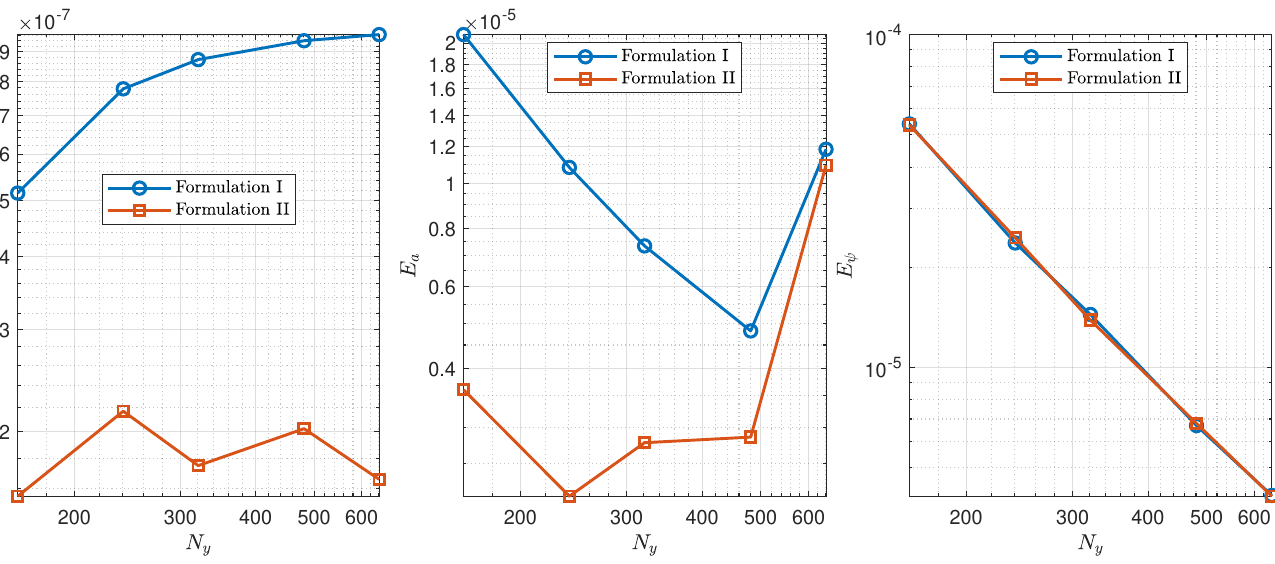}
		\caption{Extended spatial-refinement study for Formulations~I and II.
			The nonlinear benchmark uses $N_y^{\mathrm{ref}}=1281$ throughout.}
		\label{fig:final-spatial-refinement}
	\end{figure}

	The corresponding sensitivity plots are retained as numerical validation
	Figure\,\ref{fig:basic-convergence}.
	A separate higher-resolution spatial study was carried out using
	\(
	N_y=161,\;241,\;321,\;481,\;641,
	\)
	while retaining the independent $N_y^{\mathrm{ref}}=1281$ nonlinear
	benchmark.  The results are shown in
	Figure~\ref{fig:final-spatial-refinement}.
	
	The current error decreases steadily with spatial refinement.  For example,
	for Formulation~II,
	\[
	E_\psi^{II}:
	\quad
	5.35\times10^{-5}
	\;\longrightarrow\;
	2.45\times10^{-5}
	\;\longrightarrow\;
	1.39\times10^{-5}
	\;\longrightarrow\;
	6.79\times10^{-6}
	\;\longrightarrow\;
	4.10\times10^{-6}.
	\]
	Thus the current discrepancy observed on the baseline grid appears to be dominated by numerical-resolution effects at the baseline grid.
	
	The concentration errors reach a numerical floor earlier.  In particular,
	the substrate errors improve through approximately $N_y=481$ but increase
	again at $N_y=641$.  At this point, temporal interpolation, structural
	quadrature, and the finite reference solution become comparable to the
	remaining spatial error.  We therefore retain $N_y=161$ for the broad
	parameter studies and use the extended refinement only as an independent
	high-resolution validation.
	
	\subsection{Dependence on the reaction strength $\lambda_1$}
	\label{subsec:lambda1-results}
	
	To test the approximations outside the baseline regime, $\lambda_1$ was
	varied over
	\[
	\lambda_1=0.25,\;0.5,\;1,\;2,\;4,\;8,
	\]
	with all other parameters fixed.
	
\begin{figure}[htbp]
		\centering
		\includegraphics[width=0.9\linewidth, trim=0 260 0 260, clip]
		{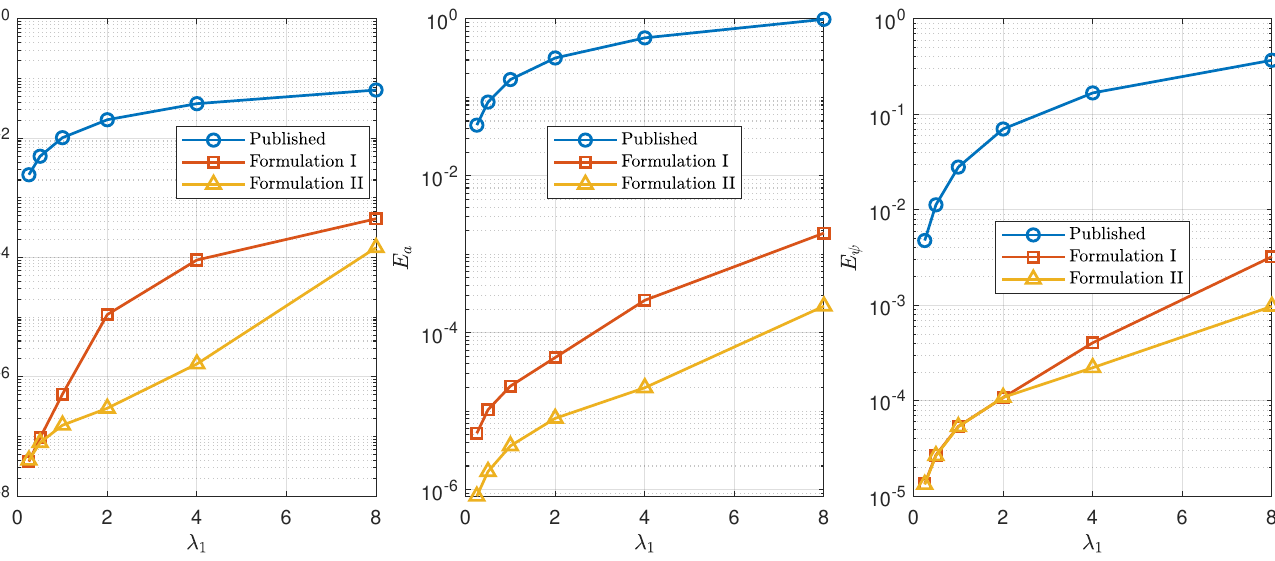}
		\caption{Dependence of the catalyst, substrate, and current errors on
			the reaction parameter $\lambda_1$.  The HPM results use $N=3$.}
		\label{fig:lambda1-final-errors}
	\end{figure}

	Figure~\ref{fig:lambda1-final-errors} shows the errors for the published
	approximation and the $N=3$ HPM approximations.  As expected, all methods
	become more challenging as the reaction rate increases.  The published
	approximation deteriorates rapidly, especially for the substrate
	concentration.  Both HPM methods remain substantially more accurate over the
	entire range.
	
	Even at $\lambda_1=8$, Formulation~II gives
	\[
	E_q^{II}=1.47\times10^{-4},
	\qquad
	E_a^{II}=2.18\times10^{-4},
	\qquad
	E_\psi^{II}=9.75\times10^{-4},
	\]
	compared with
	\[
	E_q^{I}=4.48\times10^{-4},
	\qquad
	E_a^{I}=1.86\times10^{-3},
	\qquad
	E_\psi^{I}=3.23\times10^{-3}.
	\]
	
\begin{figure}[htbp]
		\centering
		\includegraphics[width=0.9\linewidth, trim=0 260 0 260, clip]
		{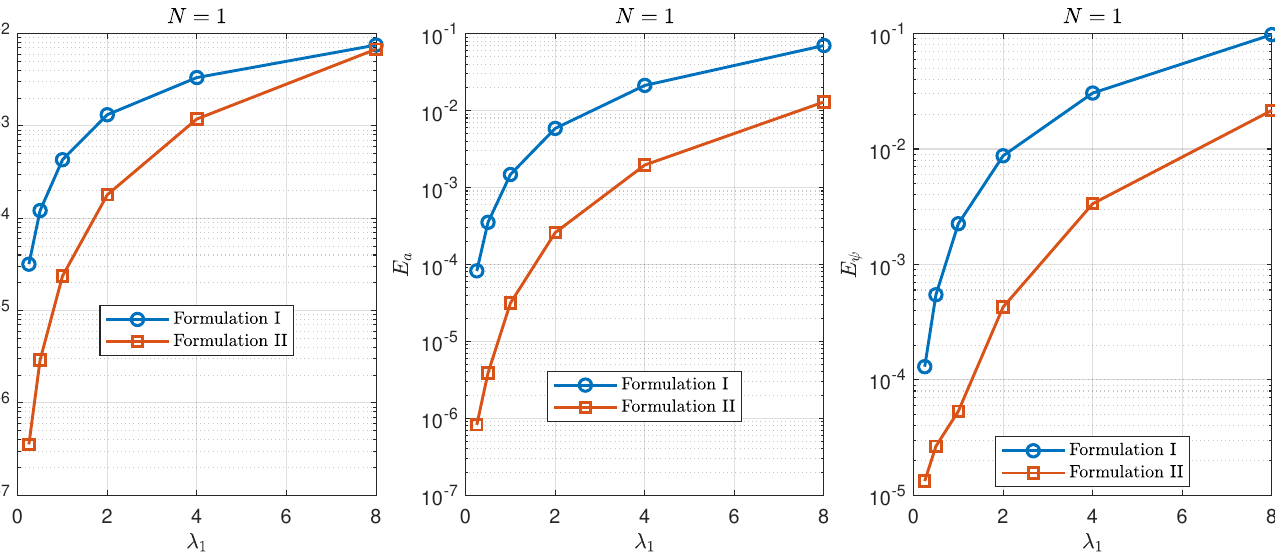}
		\caption{Low-order ($N=1$) errors as $\lambda_1$ varies.
			Formulation~II provides substantially greater accuracy than
			Formulation~I throughout the tested range.}
		\label{fig:lambda1-N1}
	\end{figure}

	The advantage of Formulation~II is especially clear at low HPM order.
	Figure~\ref{fig:lambda1-N1} compares the two approximations using only
	$N=1$.
	
	At $\lambda_1=0.25$, Formulation~II reduces the $N=1$ catalyst and substrate
	errors by factors of approximately $89.5$ and $100$, respectively.  The
	advantage decreases as the reaction becomes increasingly strong, but remains
	present at $\lambda_1=8$, where the corresponding improvement factors are
	approximately
	\(
	1.11, 5.41, 4.53
	\)
	for $q$, $a$, and $\psi$, respectively.
	The size of the next profile correction provides a useful indicator of the increasing importance of the nonlinear remainder and higher-order corrections. Over the same range,
	\[
	\frac{\|u_1\|_\infty}{\|u_0\|_\infty}
	\]
	increases from approximately
	$1.47\times10^{-4}$ at $\lambda_1=0.25$ to approximately $1.09\times10^{-1}$
	at $\lambda_1=8$.  Thus stronger reaction increases the importance of higher
	corrections, while the profile-linearized reference continues to provide a
	substantial low-order improvement.
	
	\subsection{Dependence on the denominator parameter $\rho$}
	\label{subsec:rho-results}
	
\begin{figure}[htbp]
		\centering
		\includegraphics[width=0.9\linewidth, trim=0 260 0 260, clip]
		{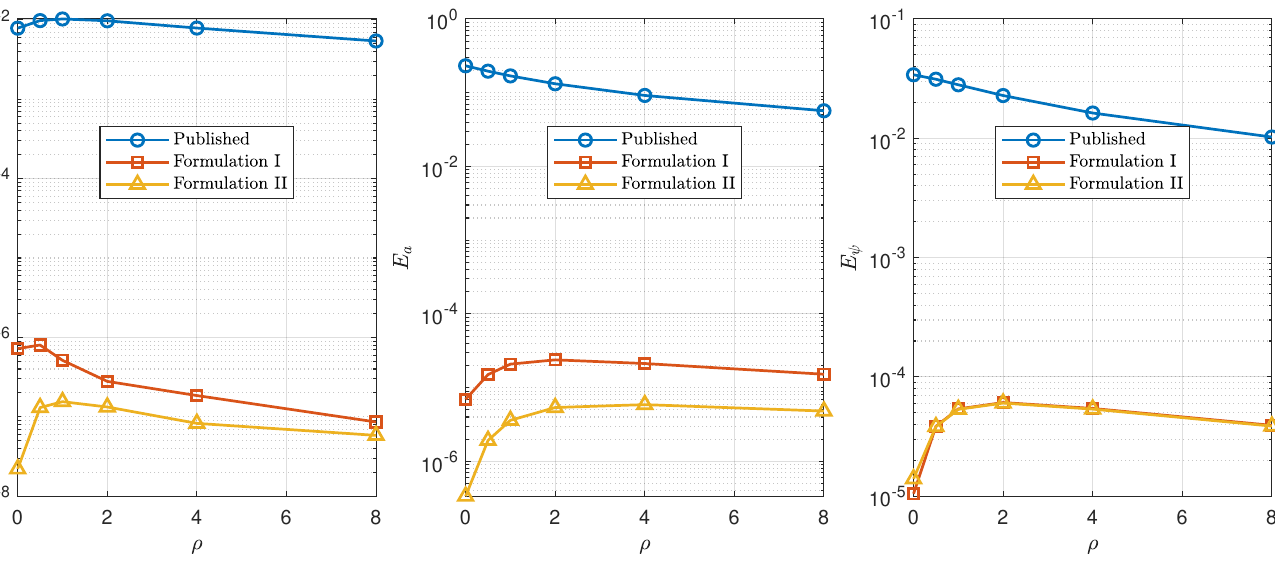}
		\caption{Dependence of the catalyst, substrate, and current errors on
			$\rho$ for the published approximation and the two $N=3$ HPM
			approximations.}
		\label{fig:rho-final-errors}
	\end{figure}

	The parameter $\rho$ was varied over
	\[
	\rho=0,\;0.5,\;1,\;2,\;4,\;8.
	\]
	Figure~\ref{fig:rho-final-errors} shows that both structure-preserving HPM
	formulations remain highly accurate throughout this range.  In contrast to
	increasing $\lambda_1$, increasing $\rho$ weakens the effective reaction and
	the published approximation improves, although its errors remain much larger
	than those of the HPM formulations.
	
\begin{figure}[htbp]
		\centering
		\includegraphics[width=0.9\linewidth, trim=0 260 0 260, clip]
		{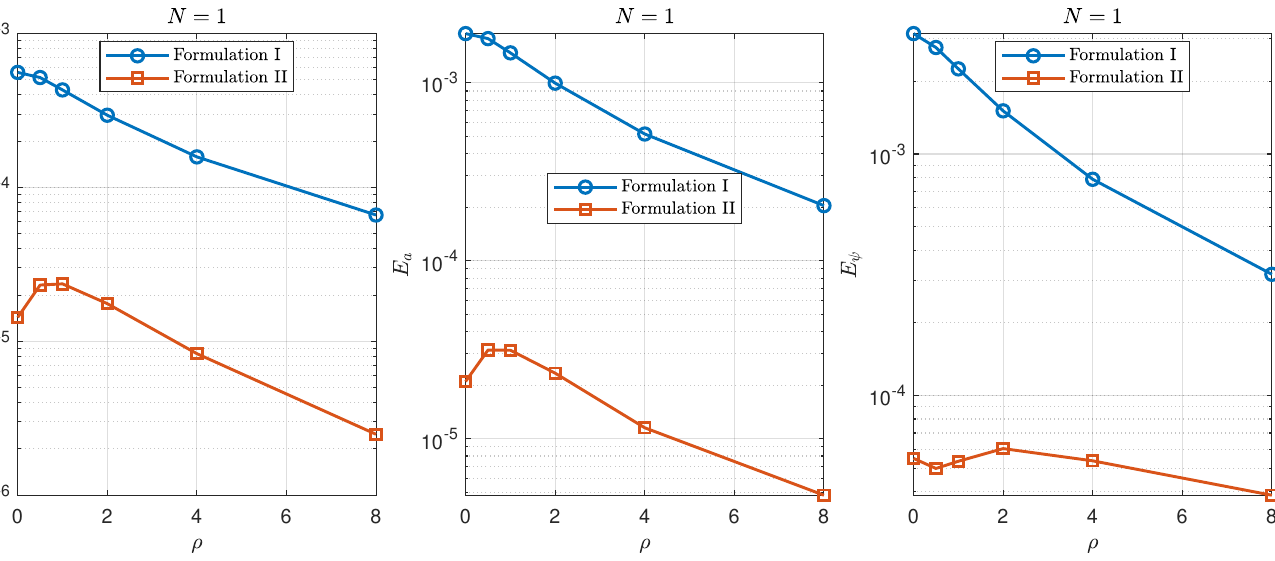}
		\caption{Low-order ($N=1$) errors as $\rho$ varies.}
		\label{fig:rho-N1}
	\end{figure}

	The low-order comparison in Figure~\ref{fig:rho-N1} again demonstrates the
	advantage of the profile-linearized formulation.
	
	Across the tested $\rho$ range, the $N=1$ Formulation~II improvement factor
	is approximately
	\(
	16.8\text{--}39.1
	\)
	for $q$,
	\(
	42.4\text{--}90.1
	\)
	for $a$, and
	\(
	8.2\text{--}57.3
	\)
	for the current.  Moreover,
	\[
	\frac{\|u_1\|_\infty}{\|u_0\|_\infty}
	\leq 2.4\times10^{-3}
	\]
	throughout the entire $\rho$ sweep. This small ratio indicates that the nonlinear correction beyond the leading profile-linearized response remains small throughout the tested \(\rho\) range.
	
	\subsection{Structural consistency}
	\label{subsec:structural-consistency}
	
	The numerical implementations were also tested against the exact structural
	relation derived earlier.  For Formulation~I, the substrate obtained by
	structural reconstruction was compared with an independent direct solution
	of the equivalent coupled frozen-coefficient MOL system.  The full-transient
	maximum difference is
	\(
	1.60\times10^{-4},
	\)
	with the largest discrepancy concentrated near the initial corner layer.
	Restricting the comparison to $\tau\geq0.05$ reduces the maximum difference
	to
	\(
	3.92\times10^{-5}.
	\)
	At the coefficient level, the post-transient discrepancies for
	Formulation~I are
	\(
	1.47\times10^{-5},
	2.46\times10^{-5},
	1.90\times10^{-7},
	3.09\times10^{-8},
	\)
	for $n=0,1,2,3$, respectively.
	
	For Formulation~II, direct computation of $v_n$ was compared with the exact
	structural reconstruction
	\(
	v_n=u_n+\mathcal H[u_n].
	\)
	The corresponding post-transient discrepancies decrease rapidly:
	\(
	2.45\times10^{-5},
	4.46\times10^{-8},
	7.93\times10^{-11}.
	\)
	These results provide numerical confirmation that both proposed algorithms
	preserve the exact difference structure to within the discretization and
	quadrature accuracy.
	
	\begin{figure}[htbp]
		\centering
		\includegraphics[width=0.9\linewidth, trim=0 260 0 260, clip]
		{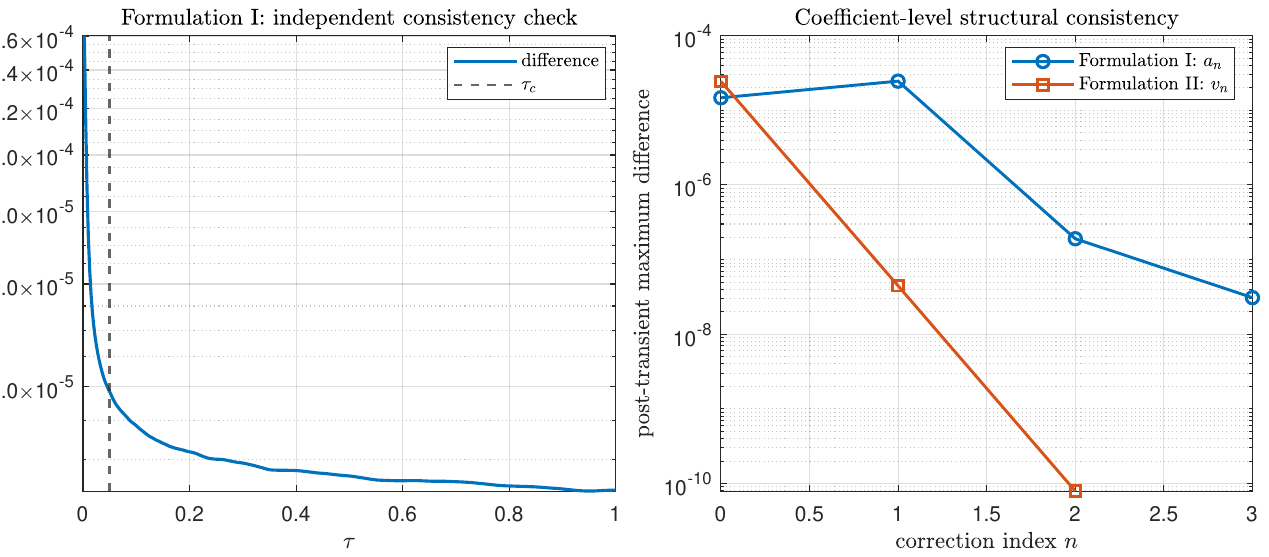}
		\caption{Structural consistency tests.  Left: difference between the
			Formulation~I structural reconstruction and the independent coupled
			MOL4 calculation.  Right: coefficient-level post-transient consistency
			for both formulations.}
		\label{fig:structural-consistency}
	\end{figure}
	
	\subsection{Computational efficiency}
	\label{subsec:efficiency}
	
	Finally, we compare the computational work required by the two low-order HPM
	approximations with that required by direct numerical solution of the
	nonlinear PDE.  Because wall-clock time is hardware dependent, two measures
	are reported:
	\begin{enumerate}\itemsep -0.1cm
		\item the median elapsed runtime from three dedicated runs; and
		\item the hardware-independent work index
		\[
		W=N_{\mathrm{state}}N_{\mathrm{RHS}},
		\]
		where $N_{\mathrm{state}}$ is the dimension of the semidiscrete ODE
		system and $N_{\mathrm{RHS}}$ is the number of right-hand-side
		evaluations performed by \texttt{ode15s}.
	\end{enumerate}
	
	The comparison is made at $N=1$, since this is the regime in which the
	profile-linearized formulation provides its most important practical
	advantage.
	
	\begin{table}[htbp]
		\centering
		\caption{Computational efficiency for the baseline problem.
			elapsed runtimes are medians of three dedicated runs.}
		\label{tab:efficiency}
		\small
		\begin{tabular}{lrrrr}
			\hline
			Method
			& $N_{\mathrm{state}}$
			& $N_{\mathrm{RHS}}$
			& $W$
			& Median runtime (s)\\
			\hline
			Formulation I, $N=1$
			& 159
			& 3352
			& $5.33\times10^{5}$
			& 2.72\\
			Formulation II, $N=1$
			& 318
			& 2522
			& $8.02\times10^{5}$
			& 4.46\\
			Nonlinear PDE
			& 2558
			& 6749
			& $1.73\times10^{7}$
			& 295.54\\
			\hline
		\end{tabular}
	\end{table}
	
	The nonlinear PDE requires approximately $32.4$ times the solver-work index
	of Formulation~I and $21.5$ times that of Formulation~II.  The corresponding
	median elapsed runtime ratios are approximately $109$ and $66$, respectively.
	The observed elapsed runtime over the three runs were approximately
	\(
	2.70\text{--}4.03\ {\rm s},
	3.91\text{--}7.06\ {\rm s},
	260\text{--}326\ {\rm s},
	\)
	for Formulation~I, Formulation~II, and the nonlinear PDE.
	
	Formulation~II is more expensive than Formulation~I at the same truncation
	order: its median elapsed runtime is approximately $1.64$ times larger and its
	solver-work index is approximately $1.50$ times larger.  However, at $N=1$
	this additional work produces error reductions of approximately
	\(
	18.2, 46.9, 42.1
	\)
	in $q$, $a$, and $\psi$, respectively.  Thus Formulation~II provides a
	particularly favorable accuracy--work tradeoff when a low-order
	approximation is desired.
	
	The kernel-based postprocessing also contributes to the HPM computational
	cost.  For the baseline calculation, Formulation~I uses $403\,200$ direct
	surface-flux quadrature samples and $400\,000$ structural-reconstruction
	samples, whereas Formulation~II uses $1\,228\,800$ surface-flux samples and
	$200\,000$ structural-reconstruction samples.  These quadrature operations
	are not included in the solver-work index $W$ and help explain the difference
	between the work-index ratios and the observed elapsed runtime ratios.
	
	\begin{figure}[htbp]
		\centering
		\includegraphics[width=0.9\linewidth, trim=0 260 0 260, clip]
		{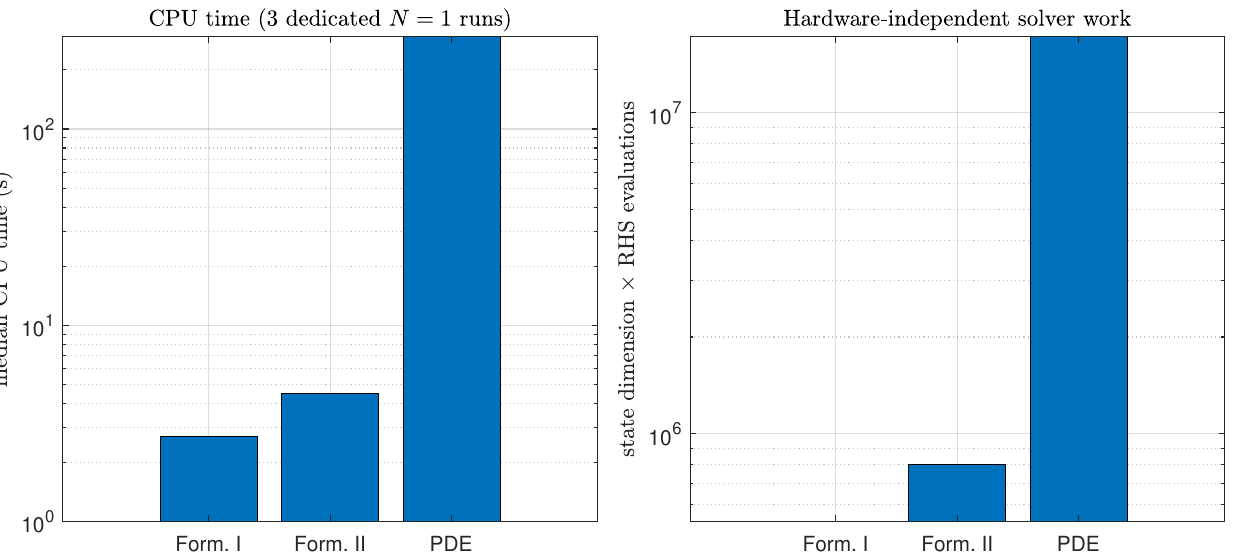}
		\caption{Computational efficiency at $N=1$.  Left: median elapsed runtime over
			three dedicated runs.  Right: hardware-independent solver-work index
			$N_{\mathrm{state}}N_{\mathrm{RHS}}$.}
		\label{fig:efficiency}
	\end{figure}
	
	\subsection{Summary of the numerical results}
	\label{subsec:numerical-summary}
	
	The numerical experiments support several conclusions.  First, the
	Green-function and fourth-order MOL realizations agree closely, providing an
	independent validation of the numerical implementation.  Second, both
	structure-preserving HPM formulations improve the published low-scan-rate
	approximation by several orders of magnitude, with the most dramatic
	improvement occurring for the substrate concentration.  Third,
	Formulation~II achieves high accuracy at substantially lower perturbation
	order than Formulation~I and remains advantageous throughout the tested
	$\lambda_1$ and $\rho$ ranges.  Fourth, the structural reconstruction checks
	confirm that the numerical algorithms preserve the exact difference
	structure derived analytically. Both HPM formulations require substantially lower semidiscrete ODE solver work and substantially less elapsed runtime than the fine-grid nonlinear PDE benchmark in the present implementation.

	\section{Discussion}
	
	The results support several conclusions about the role of structure
	preservation in the reduced electrochemical reaction--diffusion model.
	The central point is that the exact identity
	\[
	(a-q)_\tau=(a-q)_{yy}
	\]
	is independent of HPM and therefore provides a model-based consistency
	condition for any approximation of the coupled concentration fields.
	The main findings are summarized and interpreted below.
	
	\begin{enumerate}\itemsep 0.8em
		
		\item \textbf{Preserving the exact difference structure substantially
			improves the substrate approximation.}
		
		The published low-scan-rate approximation and the proposed formulations
		share the same leading catalyst field,
		\[
		q_{\mathrm P}=q_0=\bar q,
		\]
		and therefore the same leading current. The major difference appears in
		the substrate field. The published companion approximation
		\(a_{\mathrm P}\) does not generally satisfy the exact difference heat
		equation, whereas the structure-preserving zeroth-order field
		\[
		a_0=q_0+\gamma+\mathcal H[q_0]
		\]
		does. Even before higher-order HPM corrections are introduced, the
		substrate error decreases from approximately \(1.69\times10^{-1}\) for
		the published approximation to \(2.24\times10^{-2}\) for the
		structure-preserving zeroth-order field. Thus a substantial part of the
		improvement arises from enforcing the exact coupling structure rather than
		simply adding more perturbation terms.
		
		\item \textbf{Formulation~I provides a simple and highly accurate
			structure-preserving hierarchy.}
		
		Formulation~I retains a scalar constant-coefficient
		diffusion--reaction operator for every catalyst correction, while the
		corresponding substrate coefficient is reconstructed from the exact heat
		equation. This gives a common Green kernel, a direct representation of the
		surface derivative, and a straightforward recursive structure. The
		equivalent coupled frozen-coefficient formulation gives the same continuous
		approximation and therefore serves as an independent structural and
		numerical consistency check rather than as a distinct third method.
		
		\item \textbf{Formulation~II provides the largest advantage at low
			perturbation order.}
		
		The profile-linearized formulation incorporates the local sensitivity of
		the nonlinear reaction to both concentration fields into the linear
		operator. At the baseline parameters and \(N=1\), it reduces the errors in
		\(q\), \(a\), and the current relative to Formulation~I by factors of
		approximately
		\[
		18.2,\qquad 46.9,\qquad 42.1,
		\]
		respectively. By \(N=2\)--\(3\), the concentration errors approach the
		numerical error floor, so further reductions in analytical truncation error
		become increasingly difficult to distinguish from discretization,
		interpolation, quadrature, and reference-solution errors. The principal
		benefit of Formulation~II is therefore its ability to achieve high accuracy
		with only one or two correction levels.
		
		\item \textbf{The parameter dependence has a direct electrochemical
			interpretation.}
		
		From the nondimensionalization,
		\[
		\lambda_1
		=
		\frac{RT}{F}\frac{k_1c_P^0}{v},
		\qquad
		\rho=\frac{\lambda_2}{\lambda}=\frac{k_2}{k}.
		\]
		Thus \(\lambda_1\) measures the strength of the forward homogeneous
		catalytic reaction relative to the voltammetric time scale. Larger
		\(\lambda_1\) may correspond to a larger \(k_1\), a larger catalyst
		concentration \(c_P^0\), or a smaller scan rate \(v\). As
		\(\lambda_1\) increases, the higher corrections become more important and
		the low-order advantage of Formulation~II decreases, although it remains
		present over the tested range.
		
		The parameter \(\rho\), by contrast, compares the reverse homogeneous
		redox step with the irreversible removal of \(B\). Increasing
		\(\rho=k_2/k\) decreases the effective reaction term
		\[
		R(a,q)=\frac{\lambda_1aq}{1+\rho(1-q)}
		\]
		for fixed \(a\) and \(q\), producing a less demanding approximation
		regime within the reduced model. The parameter studies therefore identify electrochemical
		kinetic--transport regimes in which low-order structure-preserving
		approximations remain reliable, rather than serving only as mathematical
		sensitivity tests.
		
		\item \textbf{The current calculations and computational comparisons
			support the practical usefulness of the two formulations.}
		
		The electrode current is a more sensitive numerical quantity than the
		concentration fields because it depends on the surface derivative of
		\(q\). The differentiated Green- and heat-kernel formulas avoid direct
		numerical differentiation of the HPM concentration fields, and the spatial
		refinement study shows that the remaining current discrepancy is primarily
		a resolution effect. The structural reconstruction tests likewise confirm
		that both formulations preserve the exact difference relation to within
		the accuracy of the numerical discretization and quadrature.
		
		Both HPM formulations also require substantially less computational work
		than direct solution of the fine-grid nonlinear PDE. Formulation~I is the
		simpler and less expensive approach, whereas Formulation~II requires
		moderately more work but gives substantially greater accuracy at low
		perturbation order. The two methods therefore have complementary roles:
		Formulation~I is attractive when simplicity and low computational cost are
		priorities, while Formulation~II is preferable when high accuracy is
		required from the first one or two correction levels.
		
	\end{enumerate}
	
	Two modeling approximations should also be kept in view when interpreting
	these results. First, the analytical problem uses the frozen-potential
	boundary value \(m_0\), whereas a continuously scanned Nernstian electrode
	would generate a time-dependent surface concentration \(m(\tau)\). Second,
	the elimination of the intermediate \(b\) uses the local algebraic kinetic
	closure \eqref{eq:b-closure}. The structural identity and reconstruction
	developed here are exact for the resulting two-variable reduced model; they
	should not, however, be interpreted as exact identities of the unreduced
	four-species transport system in the absence of these approximations.
	Extension to a genuinely time-dependent Nernst boundary and to the unreduced
	transport system therefore provides a natural direction for future work.
	
	Overall, the results show that the principal benefit of the proposed
	approach comes from preserving an exact property of the reduced
	electrochemical model rather than from perturbation order alone. HPM
	provides a convenient analytical framework for exploiting this structure,
	with Formulation~I emphasizing simplicity and Formulation~II emphasizing
	high low-order accuracy.
	
	\section{Conclusions}
	
	This study identifies an exact structural property of the reduced
	reaction--diffusion model for homogeneous redox catalysis and uses it as the
	basis for constructing improved semi-analytical approximations. Because the
	catalyst and substrate equations contain the same nonlinear reaction term,
	their difference satisfies a linear heat equation exactly within the reduced
	model. This relation provides both an exact reconstruction of the substrate
	concentration from the catalyst field and a natural consistency condition
	for approximate solutions. The main contribution is therefore the
	preservation of this model structure, with HPM serving as the analytical
	framework used to construct the approximations.
	
	Two structure-preserving formulations were developed. Formulation~I retains
	a simple constant-coefficient hierarchy and is computationally inexpensive,
	whereas Formulation~II incorporates the local sensitivity of the nonlinear
	reaction through a profile-dependent linearization. Both substantially
	improve the published low-scan-rate approximation, particularly for the
	substrate concentration. The largest advantage of Formulation~II occurs at
	low perturbation order: at \(N=1\), the errors in the catalyst concentration,
	substrate concentration, and current are reduced relative to Formulation~I
	by factors of approximately \(18.2\), \(46.9\), and \(42.1\), respectively.
	Thus high accuracy can be obtained without requiring many correction levels,
	while still requiring substantially less computational work than the
	fine-grid nonlinear PDE benchmark.
	
	The parameter studies further connect the mathematical approximation to the
	underlying electrochemistry. Increasing \(\lambda_1\), corresponding to
	stronger homogeneous catalytic kinetics relative to the voltammetric time
	scale, increases the importance of higher-order corrections, whereas larger
	\(\rho=k_2/k\) weakens the effective net reaction in the reduced model and
	produces a less demanding approximation regime. These results indicate that
	structure preservation is an important criterion in the development of
	semi-analytical approximations for coupled electrochemical
	reaction--diffusion systems. Extensions to a genuinely time-dependent
	Nernst boundary condition and to the unreduced four-species transport model
	would provide natural directions for future work.

	\section*{Data availability}
	
	The MATLAB code and numerical data used to generate the results, tables,
	and figures reported in this study are publicly available through MATLAB
	Central File Exchange at
	\url{https://www.mathworks.com/matlabcentral/fileexchange/185262-structure-preserving-redox-catalysis}.
	
	\section*{CRediT authorship contribution statement}
	
	\textbf{Guangming Yao:}
	Conceptualization, Methodology, Formal analysis, Software, Validation,
	Visualization, Writing -- original draft, Writing -- review and editing,
	Supervision, Project administration.
	
	\textbf{Michael Shaw:}
	Conceptualization, Methodology, Validation, Writing -- review and editing.
	
	\textbf{Nazia Nimmy:}
	Conceptualization, Formal analysis, Investigation, Validation,
	Writing -- review and editing.
	
	\textbf{Adriana Silva:}
	Software, Investigation, Validation, Visualization,
	Writing -- review and editing.
	\section*{Funding}
	
	This research did not receive any specific grant from funding agencies in
	the public, commercial, or not-for-profit sectors.
	
	\section*{Declaration of competing interest}
	
	The authors declare that they have no known competing financial interests
	or personal relationships that could have appeared to influence the work
	reported in this paper.
	
	\section*{Declaration of generative AI and AI-assisted technologies
		in the manuscript preparation process}
	
	During the preparation of this work, the authors used ChatGPT (OpenAI) to
	assist with language editing, manuscript organization, and literature-search
	support. After using this tool, the authors reviewed and edited the content
	as needed and take full responsibility for the content of the published
	article.
	
	\bibliographystyle{unsrt}
	\bibliography{references}
	
\end{document}